\documentclass[oneside,english]{amsart}
\usepackage{amssymb,amsbsy}
\usepackage{amsthm}

\numberwithin{equation}{section} %% Comment out for sequentially-numbered
\numberwithin{figure}{section} %% Comment out for sequentially-numbered
\theoremstyle{plain}
\theoremstyle{plain}
\newtheorem{thm}{Theorem}
  \theoremstyle{plain}
  \newtheorem{prop}[thm]{Proposition}
  \theoremstyle{plain}
  \newtheorem{cor}[thm]{Corollary}
  \theoremstyle{remark}
  \newtheorem{rem}[thm]{Remark}
  \theoremstyle{plain}
  \newtheorem{lem}[thm]{Lemma}
   \theoremstyle{definition}
  \newtheorem{example}[thm]{Example}

\numberwithin{thm}{section}

\begin{document}

\title{Canonical Models For Bi-Isometries}

\author{H. Bercovici, R. G. Douglas, and C. Foias}

\subjclass{Primary: 47A45. Secondary: 47A15, 47B37}

\keywords{Bi-isometry, characteristic function, functional model, pivotal operator,
similarity}

\address{HB: Department of Mathematics\\
 Indiana University \\
Bloomington, IN 47405}

\email{bercovic@indiana.edu}

\address{RGD and CF: Department of Mathematics \\
Texas A\&M University \\
College Station, TX 77843}

\email{rdouglas@math.tamu.edu\\
foias@math.tamu.edu}

\dedicatory{We dedicate this paper to the memory Israel Gohberg,
great mathematician, wonderful human being, friend and teacher to us
all.}

\thanks{HB and RGD were supported in part by grants from the National Science
Foundation.}
\begin{abstract}
A canonical model, analogous to the one for contraction operators,
is introduced for bi-isometries, two commuting isometries on a Hilbert
space. This model involes a contractive analytic operator-valued function
on the unit disk. Various pureness conditions are considered as well
as bi-isometries for which both isometries are shifts. Several families
of examples are introduced and classified.
\end{abstract}
\maketitle

\section{Introduction}

\markboth{}{}It is difficult to overestimate the the importance of
the von Neumann-Wold theorem on the structure of isometric operators
on Hilbert space. Originally introduced in the study of symmetric
operators by von Neumann, it became the foundation for Wold's study
of stationary stochastic processes. Later, it was the starting point
for the study of contraction operators by Sz.-Nagy and the third author
as well as a key ingredient in engineering systems theory. Thus it
has had an important role in both pure mathematics and its applications.

For nearly fifty years, researchers have sought a similar structure
theory for $n$-tuples of commuting isometries [4,11,12,15,16,17,19] with varying success. In \cite{BDF} the authors rediscovered
an earlier fundamental result of Berger, Coburn and Lebow [4]
on a model for an $n$-tuple of commuting isometries and carried the
analysis beyond what the latter researchers had done. In the course
of this study, a very concrete canonical model for bi-isometries
emerged; that is for pairs of commuting isometries. This new model
is related to the canonical functional model of a contraction, but
it displays subtle differences and a new set of challenges. In this
paper we take up the systematic presentation and development of this
model.

After some preliminaries, we begin in Section 3 by examining the passage
from an $n$-isometry to an $(n+1)$-isometry showing that essentially
the main ingredient needed is a contraction in the commutant of a
completely nonunitary $n$-isometry. In the case of a bi-isometry,
this additional operator can be viewed as a contractive operator-valued
analytic function in the unit disk. It is this function that is the
heart of our canonical model. We relate the reducing subspaces of
an $n$-isometry to this construction and investigate a variety of
notions of` ``pureness'' which generalize the notion of completely
nonunitary for contractions and the results of several earlier researchers.
(See Section 3 for the details.)

In Section 4 we specialize to the case $n=1$, that is to the case of
bi-isometries, and study the extension from the first isometry to
the pair. The analytic operator function mentioned above then is the
characteristic function for the pair. Various relations between the
bi-isometry and the characteristic function are investigated. In
Section 5, this model is re-examined in the context of a functional
model; that is, one in which the abstract Hilbert spaces are
realized as Hardy spaces of vector-valued functions on the unit
disk. This representation allows one to apply techniques from
harmonic analysis in their study. In Section 6, we specialize to
bi-shifts or bi-isometries for which both isometries are shift
operators. (Note that this use of the term is not the same as that
used by earlier authors.)

In Section 7, we return to the functional model for bi-isometries
obtaining unitary invariants for them. Finally, in Section 8, several
families of bi-isometries are introduced and studied. The results
here are not exhaustive but intended to illustrate various aspects
of the earlier theory as well as the variety of possibilities presented
by bi-isometries.

At the ends of Sections 3 and 4, the connection between intertwining
maps and common invariant subspaces for bi-isometries is discussed.
This topic has already been considered in \cite{BDF-2} and further
results will be presented in another paper.

\section{Preliminaries about commuting isometries}

We will study families $\mathbb{V}=(V_{i})_{i\in I}$ of commuting
isometric operators on a complex Hilbert space $\mathfrak{H}$. A
(closed) subspace $\mathfrak{M}\subset\mathfrak{H}$ is \emph{invariant}
for $\mathbb{V}$ if $V_{i}\mathfrak{M}\subset\mathfrak{M}$ for $i\in I$;
we write $\mathbb{V}|\mathfrak{M}=(V_{i}|\mathfrak{M})_{i\in I}$
if $\mathfrak{M}$ is invariant. The invariant subspace $\mathfrak{M}$
is \emph{reducing} if $\mathfrak{M}^{\perp}$ is invariant for $\mathbb{V}$
as well. If $\mathfrak{M}$ is a reducing subspace, we have a decomposition\[
\mathbb{V}=(\mathbb{V}|\mathfrak{M})\oplus(\mathbb{V}|\mathfrak{M}^{\perp}),\]
and $\mathbb{V}|\mathfrak{M}$ is called a \emph{direct summand} of
$\mathbb{V}$. The family $\mathbb{V}$ is said to be \emph{unitary}
if each $V_{i},i\in I$, is a unitary operator. We say that $\mathbb{V}$
is \emph{completely nonunitary} or cnu if it has no unitary direct
summand acting on a space $\mathfrak{M}\ne\{0\}$. The family $\mathbb{V}$
is \emph{irreducible} if it has no reducing subspaces other than $\{0\}$
and $\mathfrak{H}$.

The following extension of the von Neumann-Wold decomposition was
proved by I. Suciu \cite{key-12}.
\begin{thm}
\label{thm:unitary+cnu(Suciu)}Let $\mathbb{V}$ be a family of commuting
isometries on $\mathfrak{H}$. There exists a unique reducing subspace
$\mathfrak{M}$ for $\mathbb{V}$ with the following properties.
\begin{enumerate}
\item [{\rm(1)}]$\mathbb{V}|\mathfrak{M}$ is unitary.
\item [{\rm(2)}]$\mathbb{V}|\mathfrak{M}^{\perp}$ is completely nonunitary.
\end{enumerate}
\end{thm}
We recall, for the reader's convenience, the construction of $\mathfrak{M}$.
We simply set\[
\mathfrak{M}=\bigcap_{N=1}^{\infty}\left[\bigcap_{k_{1},k_{2},\dots,k_{N}\in I}V_{k_{1}}V_{k_{2}}\cdots V_{k_{n}}\mathfrak{H}\right].\]
Obviously $V_{k}\mathfrak{M}\supset\mathfrak{M}$ for each $k$, and
the commutativity of $\mathbb{V}$ implies that $V_{k}\mathfrak{M}\subset\mathfrak{M}$
as well. Thus $\mathfrak{M}$ reduces each $V_{k}$ to a unitary operator.
It is then easily seen that $\mathfrak{M}$ is the largest invariant
subspace for $\mathbb{V}$ such that $\mathbb{V}|\mathfrak{M}$ is
unitary, and this immediately implies properties $(1)$ and (2), as
well as the uniqueness of $\mathfrak{M}$.

More generally, given a subset $J\subset I$, we will say that $\mathbb{V}$
is $J$-\emph{unitary} if $V_{j}$ is a unitary operator for each
$j\in J$. The family $\mathbb{V}$ is said to be $J$-\emph{pure}
if it has no $J$-unitary direct summand acting on a nonzero space.
The preceding result extends as follows.
\begin{thm}
\label{thm:unitary+cnu(given_J)}Let $\mathbb{V}=(V_{i})_{i\in I}$
be a family of commuting isometries on a Hilbert space $\mathfrak{H}$,
and let $J$ be a subset of $I$. There exists a unique reducing subspace
$\mathfrak{M}_{J}$ for $\mathbb{V}$ with the following properties.
\begin{enumerate}
\item [{\rm(1)}]$\mathbb{V}|\mathfrak{M}_{J}$ is $J$-unitary.
\item [{\rm(2)}]$\mathbb{V}|\mathfrak{M}_{J}^{\perp}$ is $J$-pure.
\end{enumerate}
\end{thm}
\begin{proof}
Let us set $\mathbb{V}_{J}=(V_{j})_{j\in J}$ and apply Theorem \ref{thm:unitary+cnu(Suciu)}
to this family. Thus we can write $\mathfrak{H}=\mathfrak{M}\oplus\mathfrak{N}$,
where $\mathfrak{M}$ is reducing for $\mathbb{V}_{J}$, $\mathbb{V}_{J}|\mathfrak{M}$
is unitary, and $\mathbb{V}_{J}|\mathfrak{N}$ is cnu. Denote by $\mathfrak{N}_{J}$
the smallest reducing subspace for $\mathbb{V}$ containing $\mathfrak{N}$,
and set $\mathfrak{M}_{J}=\mathfrak{H}\ominus\mathfrak{N}_{J}$. Since
$\mathfrak{M}_{J}$ reduces $\mathbb{V}_{J}|\mathfrak{M}$, it follows
immediately that (1) is satisfied. Moreover, if $\mathfrak{R}$ is
any reducing subspace for $\mathbb{V}$ such that $\mathbb{V}_{J}|\mathfrak{R}$
is unitary, then $\mathfrak{R}\subset\mathfrak{M}$ so that $\mathfrak{R}\perp\mathfrak{N}$
and consequently $\mathfrak{R}\perp\mathfrak{N}_{J}$ as well. We
conclude that $\mathfrak{M}_{J}$ is the largest reducing subspace
for $\mathbb{V}$ satisfying condition (1). Property (2), as well
as the uniqueness of $\mathfrak{M}_{J}$, follow from this observation.
\end{proof}
Observe that $\mathfrak{M}_{I}$ is precisely the space $\mathfrak{M}$
in Theorem \ref{thm:unitary+cnu(Suciu)}, and it is convenient to
extend our notation so that $\mathfrak{M}_{\varnothing}=\mathfrak{H}$.
We have then\[
\mathfrak{M}_{J_{1}\cup J_{2}}=\mathfrak{M}_{J_{1}}\cap\mathfrak{M}_{J_{2}},\quad J_{1},J_{2}\subset I.\]

The spaces $\mathfrak{M}_{J}$ constructed above are in fact hyperinvariant
for $\mathbb{V}$. Given two families $\mathbb{V}^{(1)}=(V_{i}^{(1)})_{i\in I}$
and $\mathbb{V}^{(2)}=(V{}_{i}^{(2)})_{i\in I}$ of commuting isometries
on $\mathfrak{H}^{(1)}$ and $\mathfrak{H}^{(2)}$, respectively,
we denote by $\mathcal{I}(\mathbb{V}^{(1)},\mathbb{V}^{(2)})$ the
collection of all bounded linear operators $X:\mathfrak{H}^{(1)}\to\mathfrak{H}^{(2)}$
satisfying the intertwining relations $XV_{i}^{(1)}=V{}_{i}^{(2)}X$
for every $i\in I$. In the special case $\mathbb{V}^{(1)}=\mathbb{V}^{(2)}=\mathbb{V}$,
we use the notation $(\mathbb{V})'=\mathcal{I}(\mathbb{V},\mathbb{V})$
for the \emph{commutant of $\mathbb{V}$.}
\begin{prop}
\label{pro:M_J-Are-hyperinvariant}Consider two families $\mathbb{V}^{(1)}=(V_{i}^{(1)})_{i\in I}$
and $\mathbb{V}^{(2)}=(V{}_{i}^{(2)})_{i\in I}$ of commuting isometries
on $\mathfrak{H}^{(1)}$ and $\mathfrak{H}^{(2)}$, respectively,
and a subset $J\subset I$. Denote by $\mathfrak{M}_{J}^{(p)}$ the
reducing subspace for $\mathbb{V}^{(p)}$ provided by the preceding
theorem for $p=1,2$. Then for every $X\in\mathcal{I}(\mathbb{V}^{(1)},\mathbb{V}^{(2)})$
we have $X\mathfrak{M}_{J}^{(1)}\subset\mathfrak{M}_{J}^{(2)}$.\end{prop}
\begin{proof}
Denote by $\mathfrak{N}_{J}^{(p)}$ the largest subspace of $\mathfrak{H}^{(p)}$
which reduces $\mathbb{V}^{(p)}$ to a unitary operator. That is,\[
\mathfrak{N}_{J}^{(p)}=\bigcap_{N=1}^{\infty}\left[\bigcap_{k_{1},k_{2},\dots,k_{N}\in J}V_{k_{1}}^{(p)}V_{k_{2}}^{(p)}\cdots V_{k_{n}}^{(p)}\mathfrak{H}^{(p)}\right],\quad p=1,2.\]
Given $X\in\mathcal{I}(\mathbb{V}^{(1)},\mathbb{V}^{(2)})$, it is
obvious from this formula that $X\mathfrak{N}_{J}^{(1)}\subset\mathfrak{N}_{J}^{(2)}$,
and therefore $X^{*}\mathfrak{N}_{J}^{(2)\perp}\subset\mathfrak{N}_{J}^{(1)\perp}$.
As noted above, $\mathfrak{M}_{J}^{(p)\perp}$ is the smallest invariant
subspace for $(V_{i}^{(p)*})_{i\in I}$ containing $\mathfrak{N}_{J}^{(p)\perp}$;
that is,\[
\mathfrak{M}_{J}^{(p)\perp}=\bigcap_{N=1}^{\infty}\left[\bigcap_{k_{1},k_{2},\dots,k_{N}\in I}V_{k_{1}}^{(p)*}V_{k_{2}}^{(p)*}\cdots V_{k_{n}}^{(p)*}\mathfrak{N}_{J}^{(p)\perp}\right],\quad p=1,2.\]
The relations $X^{*}V_{i}^{(2)*}=V_{i}^{(1)*}X^{*}$ imply now $X^{*}\mathfrak{M}_{J}^{(2)\perp}\subset\mathfrak{M}_{J}^{(1)\perp}$,
and this relation is equivalent to the conclusion of the proposition.
\end{proof}
Another useful result is the existence of a unique minimal unitary
extension for every family of commuting isometries \cite[Chapter I]{SzNFbook-1}
(see also \cite{Dou-LondMS} for a Banach space version). We review
the result briefly.
\begin{thm}
\label{thm:commutant-lifting}Let $\mathbb{V}=(V_{i})_{i\in I}$ be
a family of commuting isometries on $\mathfrak{H}$. There exists
a family $\mathbb{U}=(U_{i})_{i\in I}$ of commuting unitary operators
on a Hilbert space $\mathfrak{K}\supset\mathfrak{H}$ with the following
properties.
\begin{enumerate}
\item [{\rm(1)}]$\mathfrak{H}$ is invariant for $\mathbb{U}$ and $\mathbb{U}|\mathfrak{H}=\mathbb{V}$.
\item [{\rm(2)}]$\mathfrak{K}=\bigvee_{N=0}^{\infty}\left[\bigvee_{k_{1},k_{2},\dots,k_{N}\in I}U_{k_{1}}^{*}U_{k_{2}}^{*}\cdots U_{k_{N}}^{*}\mathfrak{H}\right].$
\end{enumerate}
If $\mathbb{U}'$ is another family of commuting unitary operators
on a space $\mathfrak{K}'\supset\mathfrak{H}$ satisfying the analogues
of conditions $(1)$ and $(2)$, then there exists a surjective isometry
$W:\mathfrak{K}\to\mathfrak{K}'$ such that $Wh=h$ for $h\in\mathfrak{H}$,
and $WU_{k}=U'_{k}W$ for $k\in I$.

\end{thm}
In equation (2) above, we use the convention that $U_{k_{1}}^{*}U_{k_{2}}^{*}\cdots U_{k_{N}}^{*}\mathfrak{H}=\mathfrak{H}$
when $N=0$.

The family $\mathbb{U}$ is called the \emph{minimal unitary extension}
of $\mathbb{V}$. In the sequel, we will denote by $\widetilde{\mathbb{V}}$
the minimal unitary extension of $\mathbb{V}$, and by $\widetilde{\mathfrak{H}}$
the space on which it acts. It is easy to verify the following \emph{commutant
extension} result. This can be deduced from the results in \cite[Chapter 1]{SzNFbook-1},
and it is proved in \cite{Dou-LondMS} for isometric operators acting
on a Banach space.
\begin{thm}
Let $\mathbb{V}^{(1)}=(V_{i}^{(1)})_{i\in I}$ and $\mathbb{V}^{(2)}=(V{}_{i}^{(2)})_{i\in I}$
be two families of commuting isometries on $\mathfrak{H}^{(1)}$ and
$\mathfrak{H}^{(2)}$, respectively, and denote by $\widetilde{\mathbb{V}^{(1)}}$
and $\widetilde{\mathbb{V}^{(2)}}$ their minimal unitary extensions.
The map $Y\mapsto X=Y|\mathfrak{H}^{(1)}$ establishes an isometric
bijection between the collection of operators $Y\in\mathcal{I}(\widetilde{\mathbb{V}^{(1)}},\widetilde{\mathbb{V}^{(2)}})$
such that $Y\mathfrak{H}^{(1)}\subset\mathfrak{H}^{(2)}$ and $\mathcal{I}(\mathbb{V}^{(1)},\mathbb{V}^{(2)})$.
\end{thm}
Indeed, given $k_{1},k_{2},\dots,k_{N}\in I$, a given operator $X\in\mathcal{I}(\mathbb{V}^{(1)},\mathbb{V}^{(2)})$
easily extends to the space $\widetilde{V_{k_{1}}^{(1)}}^{*}\widetilde{V_{k_{1}}^{(1)}}^{*}\cdots\widetilde{V_{k_{1}}^{(1)}}^{*}\mathfrak{H}^{(1)}$
by setting\[
Y\widetilde{V_{k_{1}}^{(1)}}^{*}\widetilde{V_{k_{1}}^{(1)}}^{*}\cdots\widetilde{V_{k_{1}}^{(1)}}^{*}h=\widetilde{V_{k_{1}}^{(2)}}^{*}\widetilde{V_{k_{1}}^{(2)}}^{*}\cdots\widetilde{V_{k_{1}}^{(2)}}^{*}Xh,\quad h\in\widetilde{V_{k_{1}}^{(1)}}^{*}\widetilde{V_{k_{1}}^{(1)}}^{*}\cdots\widetilde{V_{k_{1}}^{(1)}}^{*}\mathfrak{H}^{(1)},\]
and the corresponding operator $Y\in\mathcal{I}(\widetilde{\mathbb{V}},\widetilde{\mathbb{V}'})$
is obtained by taking the closure of this extension. This unique extension
of $X$ will be denoted $\widetilde{X}$. If $X$ is isometric or
unitary then so is $\widetilde{X}$. In the particular case $\mathbb{V}^{(1)}=\mathbb{V}^{(2)}=\mathbb{V}$,
the operator $X$ belongs to the commutant of $\mathbb{V}$, and its
canonical extension $\widetilde{X}\in(\widetilde{\mathbb{V}})'$ satisfies
$\widetilde{X}\mathfrak{H}\subset\mathfrak{H}$.

Irreducible families of commuting isometries have special properties.
Theorem \ref{thm:unitary+cnu(Suciu)} shows that they are either unitary
or cnu. More precisely, we have the following result.
\begin{prop}
\label{pro:irreducibles-are-special}Let $\mathbb{V}=(V_{i})_{i\in
I}$ be an irreducible family of commuting isometries on a nonzero
Hilbert space $\mathfrak{H}$. For every $i_{0}\in I$, one of the
following alternatives occurs.
\begin{enumerate}
\item [{\rm(1)}]$V_{i_{0}}$ is a scalar multiple of the identity.
\item [{\rm(2)}]$\mathbb{V}$ is $\{i_{0}\}$-pure.
\end{enumerate}
\end{prop}
\begin{proof}
Assume that (2) does not occur. Theorem \ref{thm:unitary+cnu(given_J)}
implies then that $V_{i_{0}}$ is unitary. Since the spectral projections
of $V_{i_{0}}$ reduce $\mathbb{V}$, it follows that the spectrum
of $V_{i_{0}}$ is a singleton, and therefore (1) is true.\end{proof}
\begin{prop}
\label{pro:primary_component}Let $\mathbb{V}=(V_{i})_{i\in I}$ and
$\mathbb{W}=(W_{i})_{i\in I}$ be  families of commuting isometries
on $\mathfrak{H}$ and $\mathfrak{K}$, respectively. Assume that
$\mathfrak{K}\ne\{0\}$ and $\mathbb{W}$ is irreducible. Let $(\mathfrak{M}_{\alpha})_{\alpha\in A}$
be a maximal family of pairwise orthogonal, reducing subspaces for
$\mathbb{V}$ with the property that $\mathbb{V}|\mathfrak{M}_{\alpha}$
is unitarily equivalent to $\mathbb{W}$ for every $\alpha$. If $\mathfrak{M}$
is any reducing subspace for $\mathbb{V}$ with the property that
$\mathbb{V}|\mathfrak{M}$ is unitarily equivalent to $\mathbb{W}$,
then $\mathfrak{M}\subset\bigoplus_{\alpha\in A}\mathfrak{M}_{\alpha}$.\end{prop}
\begin{proof}
This is really a general fact about representations of
$C^{*}$-algebras. We provide a proof for the sake of completeness.
Fix isometric operators $U,U_{\alpha}:\mathfrak{K}\to\mathfrak{H}$
such that
$U\mathfrak{K}=\mathfrak{M},U_{\alpha}\mathfrak{K}=\mathfrak{M}_{\alpha}$,
$U\mathbb{W}=\mathbb{V}U$ and
$U_{\alpha}\mathbb{W}=\mathbb{V}U_{\mathfrak{\alpha}}$. The
operator\[ R=\left[I-\sum_{\alpha\in
A}P_{\mathfrak{M}_{\alpha}}\right]U:\mathfrak{K}\to\mathfrak{H}\]
satisfies the relations $R\mathbb{W}=\mathbb{V}R$ and
$R\mathbb{W}^{*}=\mathbb{V}^{*}R$, and therefore $R^{*}R$ commutes
with $\mathbb{W}$ and $\mathbb{W}^{*}$. We must have then
$R^{*}R=\rho^{2}I_{\mathfrak{K}}$ for some $\rho\ge0$. If
$\rho\ne0$, then the unitary operator $U_{0}=R/\rho$ satisfies
$U_{0}\mathbb{W}=\mathbb{V}U_{0}$,
$U_{0}\mathbb{W}^{*}=\mathbb{V}^{*}U_{0}$, $U_{0}^{*}U_{\alpha}=0$,
and therefore $U_{0}\mathfrak{K}$ is a reducing space orthogonal to
each $\mathfrak{M}_{\alpha}$, contradicting the maximality of
$(\mathfrak{M}_{\alpha})_{\alpha\in A}$. Thus $\rho=0$ and the
proposition follows.
\end{proof}
It follows from this proposition that the reducing subspace\[
\mathfrak{H}_{\mathbb{W}}=\bigoplus_{\alpha\in A}\mathfrak{M}_{\alpha}\]
does not depend on the particular family $(\mathfrak{M}_{\alpha})_{\alpha\in A}$.
The restriction $\mathbb{V}|\mathfrak{H}_{\mathbb{W}}$ is an orthogonal
sum of copies of $\mathbb{W}$, while $\mathbb{V}|\mathfrak{H}_{\mathbb{W}}^{\perp}$
has no restriction to an invariant subspace that is unitarily equivalent
to $\mathbb{W}$.
\begin{prop}
\label{pro:primary_are_perp}Let $\mathbb{V}=(V_{i})_{i\in I}$,
$\mathbb{W}_{1}=(W_{i1})_{i\in I}$ and
$\mathbb{W}_{2}=(W_{i2})_{i\in I}$ be  families of commuting
isometries on $\mathfrak{H},$$\mathfrak{K}_{1}$ and
$\mathfrak{K}_{2}$, respectively. Assume that $\mathbb{W}_{1}$ and
$\mathbb{W}_{2}$ are irreducible and not unitarily equivalent. Then
the spaces $\mathfrak{\mathfrak{H}}_{\mathbb{W}_{1}}$ and
$\mathfrak{H}_{\mathbb{W}_{2}}$ are mutually orthogonal.\end{prop}
\begin{proof}
Let $\mathfrak{M}_{j}$ be a reducing subspace for $\mathbb{V}$ such
that $\mathbb{V}|\mathfrak{M}_{j}$ is unitarily equivalent to
$\mathbb{W}_{j}$ via a unitary operator
$U_{j}:\mathfrak{K}_{j}\to\mathfrak{M}_{j}$. It will suffice to show
that $\mathfrak{M}_{1}\perp\mathfrak{M}_{2}$ or, equivalently, that
the operator $R=U_{2}^{*}U_{1}$ is zero. As in the preceding proof,
irreducibility shows that $R^{*}R=\rho^{2}I_{\mathfrak{K}_{1}}$ and
$RR^{*}=\rho^{2}I_{\mathfrak{K}_{2}}$ for some constant $\rho\ge0$.
The assumption $\rho\ne0$ implies that $\mathbb{W}_{1}$ and
$\mathbb{W}_{2}$ are unitarily equivalent via the unitary operator
$R/\rho$, which is a contradiction.\end{proof}
\begin{cor}
\label{cor:primaries+remainder}Let $\mathbb{V}$ be a family of
commuting isometries on $\mathfrak{H}$, and denote by $\mathcal{F}$
a collection of mutually inequivalent irreducible families of
commuting isometries such that every irreducible direct summand of
$\mathbb{V}$ is equivalent to an element of $\mathcal{F}$. We have
\[
\mathfrak{H}=\mathfrak{H}_{0}\oplus\bigoplus_{\mathbb{W}\in\mathcal{F}}\mathfrak{H}_{\mathbb{W}},\]
where $\mathfrak{H}_{0}$ is a reducing subspace for $\mathbb{V}$
such that $\mathbb{V}|\mathfrak{H}_{0}$ has no irreducible direct
summand.
\end{cor}
When $\dim\mathfrak{H}_{0}>1$, the family
$\mathbb{V}|\mathfrak{H}_{0}$ is certainly reducible; it just cannot
be decomposed into a direct sum of irreducible families. However it
can be decomposed into a continuous direct integral of irreducibles
if $\mathfrak{H}$ is separable. A concrete example of such a
decomposition will be given in Section
\ref{sec:Examples-of-irreducible}. Direct integrals are also useful
in the proof of the following result, an early variant of which was
proved in \cite{key-12} when $I$ consists of two elements. We refer
to \cite{takesaki} for the theory of direct integrals.
\begin{prop}
\label{pro:extended-Suciu}Let $\mathbb{V}=(V_{i})_{i\in I}$ be a
finite family of commuting isometries on a Hilbert space $\mathfrak{H}$.
We can associate to each subset $J\subset I$ a reducing space $\mathfrak{L}_{J}$
for $\mathbb{V}$ with the following properties.
\begin{enumerate}
\item [{\rm(1)}]$\mathfrak{H}=\bigoplus_{J\subset I}\mathfrak{L}_{J}.$
\item [{\rm(2)}]$V_{j}|\mathfrak{L}_{J}$ is unitary for each $j\in J$.
\item [{\rm(3)}]$\mathbb{V}|\mathfrak{L}_{J}$ is $\{j\}$-pure for each
$j\notin J$.
\end{enumerate}
\end{prop}
\begin{proof}
Since $I$ is finite, $\mathfrak{H}$ can be written as an orthogonal
sum of separable reducing subspaces for $\mathbb{V}$. Thus it is
sufficient to consider the case of separable spaces $\mathfrak{H}$.
There exist a probability measure $\mu$ on $[0,1]$, a measurable
family $(\mathfrak{H}_{t})_{t\in[0,1]}$ of Hilbert spaces, and a
measurable collection $(\mathbb{V}_{t})_{t\in[0,1]}=((V_{ti})_{i\in I})_{t\in[0,1]}$
of irreducible families of commuting isometries on $\mathfrak{H}_{t}$
such that, up to unitary equivalence,\[
\mathfrak{H}=\int_{[0,1]}^{\oplus}\mathfrak{H}_{t}\, d\mu(t),\quad V_{i}=\int_{[0,1]}^{\oplus}V_{ti}\, d\mu(t),\quad i\in I.\]
Proposition \ref{pro:irreducibles-are-special} shows that for each
$t\in[0,1]$ there exists a subset $J(t)\subset I$ such that $V_{tj}$
is a scalar multiple of the identity if $j\in J(t)$, while $\mathbb{V}_{t}$
is $\{j\}$-pure for $j\notin J(t)$. It is easy to verify that the
set $\sigma_{J}=\{t\in[0,1]:J(t)=J\}$ is measurable for each $J\subset I$.
The spaces\[
\mathfrak{L}_{J}=\int_{\sigma_{J}}^{\oplus}\mathfrak{H}_{t}\, d\mu(t),\]
viewed as subspaces of $\mathfrak{H}$, satisfy the conclusion of
the proposition.
\end{proof}

\section{\label{sec:Inductive-construction-of-n+1}Inductive construction
of commuting isometries}

In this section it will be convenient to index families of commuting
isometries by ordinal numbers. Thus, given an ordinal number $n$,
an $n$-\emph{isometry} is simply a family $\mathbb{V}=(V_{i})_{0\le i<n}$
of commuting isometries on a Hilbert space.

We consider a special construction which produces an $(n+1)$-isometry
starting from an $n$-isometry $\mathbb{V}$ on $\mathfrak{H}$ and
a contraction $A\in(\mathbb{V})'$; that is, $\|A\|\le1$. Observe
that the canonical extension $\widetilde{A}\in(\widetilde{\mathbb{V}})'$
on $\widetilde{\mathfrak{H}}$ is then a contraction as well, and
therefore we can form the defect operator \[
D_{\widetilde{A}}=(I-\widetilde{A}^{*}\widetilde{A})^{1/2}\]
and the space $\mathfrak{D}=\overline{D_{\widetilde{A}}\widetilde{\mathfrak{H}}}$.
The space $\mathfrak{D}$ is reducing for $\widetilde{\mathbb{V}}$
because $D_{\widetilde{A}}$ commutes with $\widetilde{\mathbb{V}}$.
We form the space \[
\mathfrak{K}=\mathfrak{H}\oplus\mathfrak{D}\oplus\mathfrak{D}\oplus\cdots,\]
and define an $(n+1)$-isometry $\mathbb{W}_{A}=(W_{k})_{0\le k\le n}$
on $\mathfrak{K}$ as follows. For $0\le k<n$ we define\[
W_{k}=V_{k}\oplus(\widetilde{V_{k}}|\mathfrak{D})\oplus(\widetilde{V_{k}}|\mathfrak{D})\oplus\cdots,\]
while\[
W_{n}(h\oplus d_{0}\oplus d_{1}\oplus\cdots)=Ah\oplus D_{\widetilde{A}}h\oplus d_{0}\oplus d_{1}\oplus\cdots\]
if $h\in\mathfrak{H}$ and $d_{j}\in\mathfrak{D}$ for $j\in\mathbb{N}$.
It is easy to verify that $\mathbb{W}_{A}$ is in fact an $(n+1)$-isometry.
When the operator $A$ is already isometric, we have $\mathfrak{K}=\mathfrak{H}$
and $\mathbb{W}_{A}=(\mathbb{V},A)$. In this trivial sense, every
$(n+1)$-isometry is of the form $\mathbb{W}_{A}$ for some contraction
$A$ commuting with an $n$-isometry $\mathbb{V}$. We give now a
characterization of $(n+1)$-isometries which are $\{0\le k<n\}$-pure.
\begin{thm}
\label{thm:n-pure_families}Let $\mathbb{W}=(W_{k})_{0\le k\le n}$
be an $(n+1)$-isometry on $\mathfrak{K}$, where $n\ge1$. The following
conditions are equivalent.
\begin{enumerate}
\item [{\rm(1)}]$\mathbb{W}$ is $\{0\le k<n\}$-pure.
\item [{\rm(2)}]There exist a cnu $n$-isometry $\mathbb{V}$, and a contraction
$A\in(\mathbb{V})'$, such that $\mathbb{W}$ is unitarily equivalent
to $\mathbb{W}_{A}$.
\end{enumerate}
\end{thm}
\begin{proof}
Assume first that $\mathbb{W}=\mathbb{W}_{A}$, where $A$ is a contraction
in the commutant of the cnu $n$-isometry $\mathbb{V}$ on $\mathfrak{H}$.
Let $\mathfrak{M}$ be a reducing subspace for $\mathbb{W}_{A}$ with
the property that $W_{k}|\mathfrak{M}$ is unitary for all $k<n$.
Since the cnu direct summand of the $n$-isometry $(W_{k})_{0\le k<n}$
is precisely $\mathbb{V}$ viewed as acting on $\mathfrak{H}\oplus\{0\}\oplus\{0\}\oplus\cdots$,
we conclude that \[
\mathfrak{M}\subset\{0\}\oplus\mathfrak{D}\oplus\mathfrak{D}\oplus\cdots\]
and therefore $W_{n}^{*N}h\in\{0\}\oplus\mathfrak{D}\oplus\mathfrak{D}\oplus\cdots$
for every $h\in\mathfrak{M}$ and $N\ge1$. This is not possible if
$h\ne0$. Indeed, if $h=0\oplus0\oplus\cdots\oplus0\oplus d_{N}\oplus\cdots$,
and the $N$th component $d_{N}$ is the first nonzero component of
$h$, then $W_{n}^{*N}h=D_{A}d_{N}\oplus\cdots\notin\mathfrak{M}$
because $D_{A}d_{N}\ne0$.

Conversely, assume that condition (1) is satisfied. Consider the $n$-isometry
$\mathbb{W}'=(W_{k})_{0\le k<n}$, and the decomposition $\mathfrak{K}=\mathfrak{H}\oplus\mathfrak{H}^{\perp}$
into reducing subspaces for $\mathbb{W}'$ such that $\mathbb{W}'|\mathfrak{H}$
is cnu and $\mathbb{W}'|\mathfrak{H}^{\perp}$ is unitary. We denote
by $\mathbb{V}=\mathbb{W}'|\mathfrak{H}$ the cnu direct summand of
$\mathbb{W}'$, and define an operator $A$ on $\mathfrak{H}$ by
setting $A=P_{\mathfrak{H}}W_{n}|\mathfrak{H}$. Clearly $A$ is a
contraction, and the fact that $A$ commutes with $\mathbb{V}$ follows
from the fact that the unitary component $\mathfrak{H}^{\perp}$ is
obviously invariant for $W_{n}$, and therefore $A^{*}=W_{n}^{*}|\mathfrak{H}$.
Consider next the minimal unitary extension $\widetilde{\mathbb{W}'}$
which can be written as\[
\widetilde{\mathbb{W}'}=\widetilde{\mathbb{V}}\oplus(\mathbb{W}'|\mathfrak{H}^{\perp})\]
on the space $\widetilde{\mathfrak{H}}\oplus\mathfrak{H}^{\perp}$,
and the unique isometric extension $\widetilde{W_{n}}$ of $W_{n}$
in the commutant of $\widetilde{\mathbb{W}'}$. Clearly,\[
\widetilde{W_{n}}|\mathfrak{H}^{\perp}=W_{n}|\mathfrak{H}^{\perp},\]
and the compression $P_{\widetilde{\mathfrak{H}}}W_{n}|\widetilde{\mathfrak{H}}$
is precisely the contractive extension $\widetilde{A}$ of $A$ in
the commutant of $\widetilde{\mathbb{V}}$. We show next that $\widetilde{W_{n}}$
is in fact the minimal isometric dilation of $\widetilde{A}$. In
other words, the smallest invariant subspace $\mathfrak{M}$ for $\widetilde{W_{n}}$
containing $\widetilde{\mathfrak{H}}$ is $\widetilde{\mathfrak{H}}\oplus\mathfrak{H}^{\perp}$.
To prove this, observe first that, since \[
\mathfrak{M}=\bigvee_{N\ge0}\widetilde{W_{n}}^{N}\widetilde{\mathfrak{H}}\]
and $\widetilde{\mathfrak{H}}$ is invariant for $\widetilde{W_{n}}^{*}$,
the space $\mathfrak{M}$ is actually reducing for $\widetilde{W_{n}}$.
Moreover, $\widetilde{W_{i}}$ is unitary for $i<n$, and hence the
operators $\widetilde{W_{i}}^{*}$ and $\widetilde{W_{n}}$ also commute.
Thus $\mathfrak{M}$ is also a reducing space for each $\widetilde{W}_{i}$
if $i<n$. We conclude that the space $\mathfrak{M}^{\perp}\subset\mathfrak{H}^{\perp}$
reduces $\mathbb{W}$, and $\mathbb{W}'|\mathfrak{M}^{\perp}$ is
unitary. Hypothesis (1) implies that $\mathfrak{M}^{\perp}=\{0\}$.

With this preparation out of the way, we find ourselves in the familiar
territory of minimal isometric dilations \cite[Chapter II]{SzNFbook-1}.
We recall that, up to unitary equivalence, the minimal isometric dilation
of the contraction $\widetilde{A}$ is the operator $W$ defined by\[
W(h\oplus d_{0}\oplus d_{1}\oplus\cdots)=\widetilde{A}h\oplus D_{\widetilde{A}}h\oplus d_{0}\oplus d_{1}\oplus\cdots\]
on the space $\widetilde{\mathfrak{H}}\oplus\mathfrak{D}\oplus\mathfrak{D}\oplus\cdots$,
where $\mathfrak{D}=\overline{D_{\widetilde{A}}\widetilde{\mathfrak{H}}}$.
We conclude that there exists a unitary operator $U:\mathfrak{D}\oplus\mathfrak{D}\oplus\cdots\to\mathfrak{H}^{\perp}$
such that\[
(I_{\widetilde{\mathfrak{H}}}\oplus U)W=W_{n}(I_{\widetilde{\mathfrak{H}}}\oplus U).\]
The reader will verify now without difficulty that the operator $I_{\mathfrak{H}}\oplus U$
provides a unitary equivalence between $\mathbb{W}_{A}$ and $\mathbb{W}$.
\end{proof}
The preceding result shows how any $\{0\le k<n\}$-pure $(n+1)$-isometry
can be constructed from a contraction in the commutant of a cnu $n$-isometry.
General $(n+1)$-isometries are described using Theorem \ref{thm:unitary+cnu(given_J)}
with $J=\{0\le k<n\}$. We record the result below, using the lifting
concept as in \cite[Sec. II.1]{FoFrGo}.
\begin{thm}
\label{thm:general(n+1)-isoms}Let $\mathbb{W}=(W_{k})_{0\le k\le n}$
be an $(n+1)$-isometry on $\mathfrak{K}$, where $n\ge1$. There
exist reducing subspaces $\mathfrak{K}_{0}$ and $\mathfrak{K}_{1}$
for $\mathbb{W}$ with the following properties.
\begin{enumerate}
\item [(1)]$\mathfrak{K}_{0}\oplus\mathfrak{K}_{1}=\mathfrak{K}$.
\item [(2)]$W_{k}|\mathfrak{K}_{1}$ is unitary for every $k<n$.
\item [(3)]$\mathbb{W}|\mathfrak{K}_{0}$ is unitarily equivalent to $\mathbb{W}_{A}$,
where $A$ is a contraction in the commutant of a cnu $n$-isometry
$\mathbb{V}$.
\end{enumerate}
The $n$-isometry $\mathbb{V}$ on $\mathfrak{H}\subset\mathfrak{K}$
is the cnu part of $\mathbb{W}'=(W_{k})_{0\le k<n}$, and the operator
$A$ is defined by the equivalent relations\[
A=P_{\mathfrak{H}}W_{n}|\mathfrak{H},\quad A^{*}=W_{n}^{*}|\mathfrak{H}.\]
In particular, $W_{n}$ is an isometric lifting of $A$, and $\widetilde{W_{n}}$
is an isometric lifting of $\widetilde{A}$, where the extension $\widetilde{W_{n}}$
belongs to $(\widetilde{\mathbb{W}'})'$ and $\widetilde{A}\in(\widetilde{\mathbb{V}})'$.

\end{thm}
Thus, the space $\mathfrak{K}_{0}$ is simply the $\{0\le k<n\}$-pure
summand of $\mathbb{W}$.

The operators which intertwine two $(n+1)$-isometries can also be
analyzed in the context of this inductive construction. Indeed, consider
$(n+1)$-isometries $\mathbb{W}^{(p)}$ acting on $\mathfrak{K}^{(p)}$,
and the corresponding decompositions\[
\mathfrak{K}^{(p)}=\mathfrak{K}_{0}^{(p)}\oplus\mathfrak{K}_{1}^{(p)},\quad p=1,2,\]
provided by Theorem \ref{thm:general(n+1)-isoms}. In other words,
$\mathbb{W}^{(p)}|\mathfrak{K}_{0}^{(p)}$ is $\{0\le k<n\}$-pure,
and $W_{k}^{(p)}|\mathfrak{K}_{1}^{(p)}$ is unitary for $0\le k<n$.
Let us further denote by $\mathfrak{H}^{(p)}$ the cnu part of $\mathfrak{K}^{(p)}$
relative to the $n$-isometry $\mathbb{W}^{(p)\prime}=\{W_{k}^{(p)}\}_{0\le k<n}$,
and set\[
\mathbb{V}^{(p)}=\mathbb{W}^{(p)\prime}|\mathfrak{H}^{(p)},\quad A^{(p)}=P_{\mathfrak{H}^{(p)}}W_{n}^{(p)}|\mathfrak{H}^{(p)},\quad p=1,2.\]
The minimal unitary extension $\widetilde{\mathbb{W}^{(p)'}}$ of
the $n$-isometry $\mathbb{W}^{(p)'}$ acts on the space\[
\widetilde{\mathfrak{K}^{(p)}}=\widetilde{\mathfrak{K}_{0}^{(p)}}\oplus\mathfrak{K}_{1}^{(p)},\]
and we denote by $\widetilde{W_{n}^{(p)}}$ the canonical extension
of $W_{n}^{(p)}$ to this larger space. We have\[
\widetilde{W_{n}^{(p)}}=\widetilde{W_{n}^{(p)}|\mathfrak{K}_{0}^{(p)}}\oplus(W_{n}^{(p)}|\mathfrak{K}_{1}^{(p)})\]
and, as seen above, $\widetilde{W_{n}^{(p)}|\mathfrak{K}_{0}^{(p)}}$
is the minimal isometric dilation of the operator $\widetilde{A^{(p)}}$.

Any operator $X\in\mathcal{I}(\mathbb{W}^{(1)},\mathbb{W}^{(2)})$
can be represented as a matrix\[
X=\left[\begin{array}{cc}
X_{00} & X_{01}\\
X_{10} & X_{11}\end{array}\right],\]
where $X_{ij}\in\mathcal{I}(\mathbb{W}^{(1)}|\mathfrak{K}_{j}^{(1)},\mathbb{W}^{(2)}|\mathfrak{K}_{i}^{(2)})$
for $i,j\in\{0,1\}$. Theorem \ref{thm:commutant-lifting} implies
the existence of an extension $\widetilde{X}\in\mathcal{I}(\widetilde{\mathbb{W}^{(1)}},\widetilde{\mathbb{W}^{(2)}})$.
This extension will be represented by a matrix of the form\[
\widetilde{X}=\left[\begin{array}{cc}
\widetilde{X_{00}} & X_{01}\\
\widetilde{X_{10}} & X_{11}\end{array}\right]\]
relative to the decompositions $\widetilde{\mathfrak{K}^{(p)}}=\widetilde{\mathfrak{K}_{0}^{(p)}}\oplus\mathfrak{K}_{1}^{(p)}$.
\begin{prop}
\label{pro:(n+1)-intertwiners}With the above notation, the following
statements are true.
\begin{enumerate}
\item [(1)]$X_{01}=0$.
\item [(2)]The operator $Z=P_{\mathfrak{H}^{(2)}}X_{00}|\mathfrak{H}^{(1)}$
belongs to $\mathcal{I}(\mathbb{V}^{(1)},\mathbb{V}^{(2)})$ and $ZA^{(1)}=A^{(2)}Z$.
\item [(3)]The operator $B=P_{\widetilde{\mathfrak{H}_{0}^{(2)}}}\widetilde{X}$
belongs to $\mathcal{I}(\widetilde{\mathbb{V}^{(1)}},\widetilde{\mathbb{V}^{(2)}})$,
$B\widetilde{W_{n}^{(1)}}=\widetilde{A^{(2)}}B$, and $B\mathfrak{H}^{(1)\perp}=\{0\}$.
\end{enumerate}
\end{prop}
\begin{proof}
Part (1) follows from Proposition \ref{pro:M_J-Are-hyperinvariant}
applied to the particular case $J=\{0\le k<n\}$. The intertwining
properties of $Z$ in part (2) follow from the fact that the space
$\mathfrak{H}^{(p)}$ is reducing for $\mathbb{W}^{(p)\prime}$ and
invariant for $W_{n}^{*}$. In other words, we can use the fact that,
relative to the decompositions $\mathfrak{K}^{(p)}=\mathfrak{H}^{(p)}\oplus\mathfrak{H}^{(p)\perp}$,
the relevant operators have matrices of the form\[
X=\left[\begin{array}{cc}
Z & 0\\*
 & *\end{array}\right],\ W_{n}^{(p)}=\left[\begin{array}{cc}
A^{(p)} & 0\\*
 & *\end{array}\right],\ W_{k}^{(p)}=\left[\begin{array}{cc}
V_{k}^{(p)} & 0\\
0 & *\end{array}\right],\quad0\le k<n.\]
For part (3) we may assume that $\mathbb{W}^{(p)},p=1,2,$ are $\{0\le k<n\}$-pure.
Hence part (3) follows from similar considerations.
\end{proof}
In the lifting framework of \cite[Sec. II.1]{FoFrGo}, the operator
$\widetilde{X}$ is said to be a \emph{lifting} of $B$, and this
lifting is contractive if $\|\widetilde{X}\|\le1$. A natural question
arises: given a contraction $B$ satisfying the requirements of Proposition
\ref{pro:(n+1)-intertwiners}(3), can one construct a contractive
lifting $\widetilde{X}\in\mathcal{I}(\widetilde{\mathbb{W}^{(1)}},\widetilde{\mathbb{W}^{(2)}})$?
If one pursues the more modest goal of finding a contractive lifting
$\widetilde{X}\in\mathcal{I}(\widetilde{\mathbb{W}_{n}^{(1)}},\widetilde{\mathbb{W}_{n}^{(2)}})$,
the answer is in the affirmative, and a parametrization of all such
contractive liftings can be extracted from \cite[Chapter VI]{FoFrGo}.
We describe the result below, under the additional assumption that
$\mathbb{W}^{(2)}$ is $\{0\le k<n\}$-pure. In the notation adopted
in this section, this amounts to the requirement that $\mathfrak{K}_{1}^{(2)}=\{0\}$.
\begin{prop}
\label{pro:parametrization}With the preceding notation, assume that
$B\in\mathcal{I}(\widetilde{W_{0}^{(1)}},A^{(2)})$ is an operator
of norm $\le1$. The set of contractive liftinge $\widetilde{X}\in\mathcal{I}(\widetilde{W_{1}^{(1)}},\widetilde{W_{1}^{(2)}})$
of $B$ is parametrized by $($that is, it is in a \emph{canonical} bijection
with$)$ the set of all contractive analytic functions $R:\mathbb{D}\to\mathcal{L}(\mathfrak{G},\mathfrak{G}')$,
where the spaces $\mathfrak{G}$ and $\mathfrak{G}'$ are given by
the formulas\[
\mathfrak{G}=\overline{D_{B}\widetilde{\mathfrak{H}^{(1)}}}\ominus\overline{D_{B}\widetilde{W_{1}^{(1)}}\widetilde{\mathfrak{H}^{(1)}}},\]
\begin{eqnarray*}
\mathfrak{G}' & = & \left[\overline{(\widetilde{W_{1}^{(2)}}-A^{(2)})B\widetilde{\mathfrak{H}^{(1)}}}\oplus\overline{D_{B}\widetilde{\mathfrak{H}^{(1)}}}\right]\\
 &  & \ominus\overline{\{(\widetilde{W_{1}^{(2)}}-A^{(2)})B\widetilde{h^{(1)}}\oplus D_{B}\widetilde{h^{(1)}}:\widetilde{h^{(1)}}\in\widetilde{\mathfrak{H}^{(1)}}\}},\end{eqnarray*}
and where $D_{B}=(I-B^{*}B)^{1/2}$.
\end{prop}
One of the liftings considered above will yield an operator $X\in\mathcal{I}(\mathbb{W}^{(1)},\mathbb{W}^{(2)})$
only when it also satisfies the conditions $\widetilde{X}\widetilde{W_{k}^{(1)}}=\widetilde{W_{k}^{(2)}}\widetilde{X}$
for $0\le k<n$, and $B$ itself is subject to the supplementary conditions\[
B\mathfrak{H}^{(1)\perp}=\{0\},\ B\mathfrak{H}^{(1)}\subset\mathfrak{H}^{(2)}\ B\in\mathcal{I}(\widetilde{\mathbb{V}^{(1)}},\widetilde{\mathbb{V}^{(2)}}).\]

We continue the discussion now under the assumption that the
operator $B$ does satisfy these additional conditions. The fact that
$B\in\mathcal{I}(\widetilde{\mathbb{V}^{(1)}},\widetilde{\mathbb{V}^{(2)}})$
is easily seen to imply that
$D_{B}\in\mathcal{I}(\widetilde{\mathbb{V}^{(1)}},\widetilde{\mathbb{V}^{(2)}})$.
Using the notation in Proposition \ref{pro:parametrization}, these
intertwining conditions imply\begin{equation}
\widetilde{V_{k}^{(1)}}\mathfrak{G}\subset\mathfrak{G}\text{ and
}(\widetilde{V_{k}^{(2)}}\oplus\widetilde{V_{k}^{(1)}})\mathfrak{G}'\subset\mathfrak{G}'\text{
for }0\le k<n.\label{eq:G-G-invariance}\end{equation} Some
additional application of techniques from of \cite[Chapter
VI]{FoFrGo} yields the following result.
\begin{prop}
\label{pro:B1}With the above notation, assume that $\mathbb{W}^{(2)}$
is $\{0\le k<n\}$-pure. The set of contractions in $\mathcal{I}(\mathbb{W}^{(1)},\mathbb{W}^{(2)})$
can be parametrized by pairs $(B,R)$, where $B\in\mathcal{I}(\widetilde{W_{n}^{(1)}},\widetilde{A^{(2)}})$
is a contraction satisfying the conditions in Proposition \emph{\ref{pro:(n+1)-intertwiners}(3)}
and $R$ is a \emph{parameter} as in Proposition \emph{\ref{pro:parametrization}}
satisfying the additional conditions\[
(\widetilde{V_{k}^{(1)}}\oplus\widetilde{V_{k}^{(2)}})R(z)=R(z)\widetilde{V_{k}^{(1)}}|\mathfrak{G},\quad0\le k<n,z\in\mathbb{D}.\]

\end{prop}
These results enable one to begin a systematic study of the
invariant subspaces of bi-isometries. This study was already started
in [3] and it will be continued in a forthcoming paper.
\begin{rem}
We emphasize again that the preceding result does not require that
$\mathbb{W}^{(1)}$ is $\{0\le k<n\}$ is pure.
\end{rem}

\section{The structure of bi-isometries}

For the remainder of this paper, we focus on bi-isometries $\mathbb{W}=(W_{0},W_{1})$
on a Hilbert space $\mathfrak{K}$. In view of Theorem \ref{thm:n-pure_families},
Theorem \ref{thm:unitary+cnu(given_J)} takes the following form when
$J=\{0\}$.
\begin{prop}
\label{pro:first_decomposition}Consider a bi-isometry
$\mathbb{W}=(W_{0},W_{1})$ on $\mathfrak{K}$, let
$\mathfrak{K}=\mathfrak{H}\oplus\mathfrak{H}^{\perp}$ be the von
Neumann-Wold decomposition relative to $W_{0}$, so that
$V_{0}=W_{0}|\mathfrak{H}$ is a unilateral shift and
$W_{0}|\mathfrak{H}^{\perp}$ is unitary. Denote by
$\widetilde{W_{0}}=\widetilde{V_{0}}\oplus(W_{0}|\mathfrak{H}^{\perp})\in\mathcal{L}(\widetilde{\mathfrak{K}})=\mathcal{L}(\widetilde{\mathfrak{H}}\oplus\mathfrak{H}^{\perp})$
the minimal unitary extension of $W_{0}$, and denote by
$\widetilde{W_{1}}\in\mathcal{L}(\widetilde{\mathfrak{K}})$ the
unique isometric extension of $W_{1}$ which commutes with
$\widetilde{W_{0}}$. Define\[
\mathfrak{M}=\bigvee_{k=0}^{\infty}\widetilde{W_{1}}^{k}\widetilde{\mathfrak{H}},\quad\mathfrak{N}=\widetilde{\mathfrak{K}}\ominus\mathfrak{M}.\]
Then the subspace $\mathfrak{N}\subset\mathfrak{H}^{\perp}$ is
reducing for $\mathbb{W}$, and $W_{0}|\mathfrak{N}$ is unitary.
Moreover, $\mathfrak{N}$ is the largest reducing subspace for
$\mathbb{W}$ with the property that $W_{0}|\mathfrak{N}$ is
unitary.\end{prop}
\begin{cor}
\label{cor:minimal_coiso_ext}With the notation of the preceding result,
the following assertions are equivalent.
\begin{enumerate}
\item [{\rm(1)}]$\mathfrak{N}=\{0\}$.
\item [{\rm(2)}]$\mathbb{W}$ is $\{0\}$-pure
\item [{\rm(3)}]The operator $\widetilde{W_{1}}^{*}$ is the minimal coisometric
extension of $\widetilde{W_{1}}^{*}|\widetilde{\mathfrak{H}}$.
\end{enumerate}
\end{cor}
The particular case of Proposition \ref{pro:extended-Suciu} for bi-isometries
can be proved by repeated application of Proposition \ref{pro:first_decomposition}.
This result was obtained first in the case of doubly commuting isometries
in \cite{key-11}; the general case appears in \cite{Gasp} (see also
\cite{DPOpov3} for another proof).
\begin{cor}
\label{cor:four_space_decomposition}Consider a bi-isometry $\mathbb{W}=(W_{0},W_{1})$
on $\mathfrak{K}$. There exist unique reducing subspaces $\mathfrak{K}_{00},\mathfrak{K}_{11},\mathfrak{K}_{01},\mathfrak{K}_{11}$
for $\mathbb{W}$ with the following properties.
\begin{enumerate}
\item [{\rm(1)}]$W_{0}|\mathfrak{K}_{01}$ is a shift and $W_{1}|\mathfrak{K}_{01}$
is unitary.
\item [{\rm(2)}]$W_{0}|\mathfrak{K}_{10}$ is unitary and $W_{1}|\mathfrak{K}_{10}$
is a shift.
\item [{\rm(3)}]$W_{0}|\mathfrak{K}_{11}$ and $W_{1}|\mathfrak{K}_{11}$
are unitary.
\item [{\rm(4)}]There is no nonzero reducing subspace $\mathfrak{N}\subset\mathfrak{K}_{00}$
for $\mathbb{W}$ such that either $W_{0}|\mathfrak{N}$ or $W_{1}|\mathfrak{N}$
is unitary.
\item [{\rm(5)}]$\mathfrak{K}=\mathfrak{K}_{00}\oplus\mathfrak{K}_{01}\oplus\mathfrak{K}_{10}\oplus\mathfrak{K}_{11}$.
\end{enumerate}
\end{cor}
\begin{proof}
Proposition \ref{pro:first_decomposition} yields a decomposition
$\mathfrak{K}=\mathfrak{N}^{\perp}\oplus\mathfrak{N}$ into reducing
subspaces for $\mathbb{W}$ such that $W_{0}|\mathfrak{N}$ is unitary
and there is no reducing subspace
$\mathfrak{N}'\subset\mathfrak{N}^{\perp}$ for $\mathbb{W}$ such
that $W_{0}|\mathfrak{N}'$ is unitary. Apply this result with the
pair $\mathbb{W}$ replaced by
$(W_{1}|\mathfrak{N},W_{0}|\mathfrak{N})$ and
$(W_{1}|\mathfrak{N}^{\perp},W_{0}|\mathfrak{N}^{\perp})$,
respectively, to obtain decompositions
$\mathfrak{N}=\mathfrak{K}_{10}\oplus\mathfrak{K}_{11}$ and
$\mathfrak{N}^{\perp}=\mathfrak{K}_{00}\oplus\mathfrak{K}_{01}$,
respectively, into sums of reducing subspaces such that
$W_{1}|\mathfrak{K}_{11}$ and $W_{1}|\mathfrak{K}_{01}$ are unitary.
Moreover, there is no nontrivial reducing subspace $\mathfrak{M}$
for $\mathbb{W}$ contained in either $\mathfrak{K}_{10}$ or
$\mathfrak{K}_{00}$ such that $W_{1}|\mathfrak{M}$ is unitary. We
leave the remaining verifications to the interested reader.
\end{proof}
Consider a bi-isometry $\mathbb{W}=(W_{0},W_{1})$ on the Hilbert
space $\mathfrak{K}$. As in Proposition \ref{pro:first_decomposition},
we consider the Wold decomposition $\mathfrak{K}=\mathfrak{H}\oplus\mathfrak{H}^{\perp}$
for $W_{0}$, with\[
\mathfrak{H}_{0}=\bigoplus_{k=0}^{\infty}W_{0}^{k}\mathfrak{E},\quad\mathfrak{E}=\ker W_{0}^{*}=\mathfrak{K}\ominus W_{0}\mathfrak{K},\]
 and we set $V_{0}=W_{0}|\mathfrak{H}$ and $A=P_{\mathfrak{H}}W_{1}|\mathfrak{H}$.
Thus, $V_{0}$ is a unilateral shift and, as observed earlier, $A$
is a contraction in the commutant of $V_{0}$. We will call $(V_{0},A)$
the \emph{characteristic pair} associated to the bi-isometry $\mathbb{W}$.
Thus, the characteristic pair is simply formed by a unilateral shift
and a contraction in its commutant. The concept of unitary equivalence
for these objects is natural: two such pairs are said to be unitarily
equivalent if they are conjugated by a unitary operator (the same
for the two operators of the pair).

The pair $(W_{1},W_{0})$ is also a bi-isometry, and the above procedure
associates to it a characteristic pair. The characteristic pairs of
$(W_{0},W_{1})$ and $(W_{1},W_{0})$ are not unitarily equivalent
in general.

For future reference, we restate Theorem \ref{thm:n-pure_families}
for the special case $n=1$, that is, the case of bi-isometries.
\begin{prop}
\label{pro:H+D+D+D_model}Let $V_{0}\in\mathcal{L}(\mathfrak{H})$ be
a unilateral shift, and $A\in\{V_{0}\}'$ a contraction. Denote by
$\widetilde{V_{0}}\in\mathcal{L}(\widetilde{\mathfrak{H}_{0}})$ the
minimal unitary extension of $V_{0}$, let
$\widetilde{A}\in\{\widetilde{V_{0}}\}'$ be the extension of
$A$, and set $D=(I-\widetilde{A}^{*}\widetilde{A})^{1/2}$,
$\mathfrak{D}=(D\widetilde{\mathfrak{H}_{0}})^{-}$.
\begin{enumerate}
\item [{\rm(1)}]The space $\mathfrak{D}$ is reducing for $\widetilde{V_{0}}$.
\item [{\rm(2)}]Define the Hilbert space \[
\mathfrak{K}=\mathfrak{H}\oplus\mathfrak{D}\oplus\mathfrak{D}\oplus\cdots,\]
and the operators $W_{0},W_{1}\in\mathcal{L}(\mathfrak{K})$ by \[
W_{0}(h\oplus d_{0}\oplus
d_{1}\oplus\cdots)=V_{0}h\oplus\widetilde{V_{0}}d_{0}\oplus\widetilde{V_{0}}d_{1}\oplus\cdots,\]
\[
W_{1}(h\oplus d_{0}\oplus d_{1}\oplus\cdots)=Ah\oplus Dh\oplus d_{0}\oplus d_{1}\oplus\cdots.\]
Then $(W_{0},W_{1})$ is a $\{0\}$-pure bi-isometry whose characteristic
pair is unitarily equivalent to $(V_{0},A)$.
\end{enumerate}
\end{prop}
We collect in the following statement some basic properties of the
characteristic pair. These follow immediately from the results in
Section \ref{sec:Inductive-construction-of-n+1}.
\begin{prop}
\label{pro:model_data} Let $\mathbb{W}=(W_{0},W_{1})$ and $\mathbb{W}'=(W_{0}',W_{1}')$
be two bi-isometries with characteristic pairs $(V_{0},A)$ and $(V_{0}',A')$,
respectively.
\begin{enumerate}
\item [{\rm(1)}]The characteristic pair of $\mathbb{W}\oplus\mathbb{W}'$
is $(V_{0}\oplus V'_{0},A\oplus A')$.
\item [{\rm(2)}]If $\mathbb{W}$ is unitarily equivalent to $\mathbb{W}',$
then $(V_{0},A)$ is unitarily equivalent to $(V_{0}',A')$.
\item [{\rm(3)}]Assume in addition that $\mathbb{W}$ and $\mathbb{W}'$
are $\{0\}$-pure. If $(V_{0},A)$ is unitarily equivalent to $(V_{0}',A')$,
then $\mathbb{W}$ is unitarily equivalent to $\mathbb{W}'$.
\item [{\rm(4)}]For every pair $(V_{0},A)$, where $V_{0}$ is a unilateral
shift and $A\in\{V_{0}\}'$ is a contraction, there exists a bi-isometry
$\mathbb{W}$ such that $(V_{0},A)$ is the characteristic pair associated
to $\mathbb{W}$. This bi-isometry can be chosen to be $\{0\}$-pure.
\end{enumerate}
\end{prop}
Part (2) of this proposition characterizes the reducing subspaces
of a $\{0\}$-pure bi-isometry in terms of its characteristic pair.
General invariant subspaces of a bi-isometry are not characterized
as easily. One difficulty is the fact that the restriction of a $\{0\}$-pure
bi-isometry to an invariant subspace is not always $\{0\}$-pure.
Assume then that we start with a $\{0\}$-pure bi-isometry $\mathbb{W}$
on $\mathfrak{K}$, $\mathfrak{K}'\subset\mathfrak{K}$ is an invariant
subspace for $\mathbb{W}$, and $\mathbb{W}'=\mathbb{W}|\mathfrak{K}'$.
The inclusion operator $X\in\mathcal{L}(\mathfrak{K}',\mathfrak{K})$
is obviously an isometry in $\mathcal{I}(\mathbb{W}',\mathbb{W})$.
Conversely, given an isometric intertwining between bi-isometries
$X\in\mathcal{I}(\mathbb{W}^{(1)},\mathbb{W})$, the range of $X$
is an invariant subspace for $\mathbb{W}$. Thus the description of
invariant subspaces for bi-isometries can be achieved by understanding
the structure of isometric operators intertwining two bi-isometries.
In the terminology of Proposition \ref{pro:B1}, one needs to find
the \emph{parameters $R$ }which give rise to isometric liftings $\widetilde{X}$
of a given contraction $B$. We presented in \cite{BDF-2} some general
results concerning this problem, and further results will appear in
a forthcoming paper.

\section{Functional representation}

The data in a characteristic pair $(V_{0},A)$ on $\mathfrak{H}$ can
alternately be encoded in a contractive analytic operator-valued
function on the unit disk $\mathbb{D}$. Set
$\mathfrak{E}=\mathfrak{H}\ominus V_{0}\mathfrak{H}$, and define
operators $\Theta_{k}\in\mathcal{L}(\mathfrak{E})$ as follows:\[
\Theta_{k}=P_{\mathfrak{E}}V_{0}^{*k}A|\mathfrak{E},\quad k\ge0.\]
We can then associate to the pair $(V_{0},A)$ the operator-valued
analytic function\[
\Theta(z)=\sum_{k=0}^{\infty}z^{k}\Theta_{k}=P_{\mathfrak{E}}(I-zV_{0}^{*})^{-1}A|\mathfrak{E},\quad|z|<1.\]
When $(V_{0},A)$ is the characteristic pair of a bi-isometry
$\mathbb{W}$, $\Theta$ will be called the \emph{characteristic
function} of $\mathbb{W}$; we will use the notation
$\Theta=\Theta_{\mathbb{W}}$ when it is necessary. If $\Theta$ is
the characteristic function of $\mathbb{W}=(W_{0},W_{1})$, then its
coefficients satisfy\[
\Theta_{k}=P_{\mathfrak{E}}V_{0}^{*k}A|\mathfrak{E}=P_{\mathfrak{E}}W_{0}^{*k}P_{\mathfrak{H}_{0}}W_{1}|\mathfrak{E}=P_{\mathfrak{E}}W_{0}^{*k}W_{1}|\mathfrak{E},\quad
k\ge0,\] since
$P_{\mathfrak{E}}W_{0}^{*k}P_{\mathfrak{H}_{0}}=P_{\mathfrak{E}}W_{0}^{*k}$.
In particular, the constant coefficient \[
\Theta_{0}=\Theta(0)=P_{\mathfrak{E}}W_{1}|\mathfrak{E}=(W_{1}^{*}|\mathfrak{E})^{*}\]
is precisely the pivotal operator associated with the pair
$(W_{1},W_{0})$, as defined in \cite{BDF}. For the convenience of
the reader, we recall that the adjoint of the \emph{pivotal
operator} associated with a bi-isometry $(W_{0},W_{1})$ is defined
as $W_{0}^{*}|\ker(W_{1}^{*})$.

The operator $\Theta(z)$ is a contraction for $|z|<1$; in fact\[
\sup_{|z|<1}\|\Theta(z)\|=\|A\|,\]
where $A^{*}=W_{1}^{*}|\mathfrak{H}$. Unitary equivalence of characteristic
functions is defined in the natural way: $\Theta$ is unitarily equivalent
to $\Theta'$ if $U\Theta(z)=\Theta'(z)U,z\in\mathbb{D},$ for some
unitary operator $U$. It may be useful to contrast this notion of
unitary equivalence with the weaker notion of \emph{coincidence,}
which is the appropriate concept in the study of functional models
for contractions \cite{SzNFbook-1}. Two operator-valued analytic
functions $\Theta$ and $\Theta'$ are said to coincide if there exist
unitary operators $U$ and $V$ such that $U\Theta(z)=\Theta'(z)V$
for all $z\in\mathbb{D}$.

Proposition \ref{pro:model_data} can now be reformulated as follows.
\begin{cor}
\label{cor:char_func}Let $\mathbb{W}$ and $\mathbb{W}'$ be two
bi-isometries with characteristic functions $\Theta$ and $\Theta'$,
respectively.
\begin{enumerate}
\item [{\rm(1)}]The characteristic function of $\mathbb{W}\oplus\mathbb{W}'$
is given by $\Theta(z)\oplus\Theta'(z)$ for $z\in\mathbb{D}$.
\item [{\rm(2)}]If $\mathbb{W}$ is unitarily equivalent to $\mathbb{W}'$,
then $\Theta$ is unitarily equivalent to $\Theta'$.
\item [{\rm(3)}]Assume in addition that $\mathbb{W}$ and $\mathbb{W}'$
are $\{0\}$-pure. If $\Theta$ is unitarily equivalent to $\Theta'$
then $\mathbb{W}$ is unitarily equivalent to $\mathbb{W}'$.
\item [{\rm(4)}]For every contractive analytic function $\Theta:\mathbb{D}\to\mathcal{L}(\mathfrak{E})$,
there exists a $\{0\}$-pure bi-isometry $\mathbb{W}$ such that $\Theta_{\mathbb{W}}$
is unitarily equivalent to $\Theta$.
\end{enumerate}
\end{cor}
In order to translate the result of Proposition
\ref{pro:H+D+D+D_model} into function theoretical terms we need some
notation. First, given a separable, complex Hilbert space
$\mathfrak{E}$, we denote as usual by $H^{2}(\mathfrak{E})$ the
Hilbert space of all square summable power series with coefficients
in $\mathfrak{E}$. Given a contractive analytic function
$\Theta:\mathbb{D}\to\mathcal{L}(\mathfrak{E})$, the analytic
Toeplitz operator $T_{\Theta}\in\mathcal{L}(H^{2}(\mathfrak{E}))$ is
defined simply as pointwise multiplication by $\Theta$. The
particular case $\Theta(z)=zI_{\mathfrak{E}}$ yields the unilateral
shift $S_{\mathfrak{E}}$. The minimal unitary extension of
$S_{\mathfrak{E}}$ is the bilateral shift $U_{\mathfrak{E}}$ on the
Hilbert space $L^{2}(\mathfrak{E})$ of all square summable Laurent
series with coefficients in $\mathfrak{E}$. The extension of
$T_{\Theta}$ which commutes with $U_{\mathfrak{E}}$ is the Laurent
operator $L_{\Theta}$ with symbol $\Theta$.

Now, the space $L^{2}(\mathfrak{E})$ can also be viewed as the space
of square integrable $\mathfrak{E}$-valued functions
$f:\mathbb{T}=\partial\mathbb{D}\to\mathfrak{E}$. When viewed in
this manner, the operator $L_{\Theta}$ is given by\[
(L_{\Theta}f)(\zeta)=\Theta(\zeta)f(\zeta)\] for almost every
$\zeta\in\mathbb{T}$, where the strong operator limit\[
\Theta(\zeta)=\lim_{r\uparrow1}\Theta(r\zeta)\] exists almost
everywhere. Similarly, the operator
$D=(I-L_{\Theta}^{*}L_{\Theta})^{1/2}$ is given as a multiplication
operator by the strongly measurable operator-valued function\[
\Delta(\zeta)=(I-\Theta(\zeta)^{*}\Theta(\zeta))^{1/2},\quad\zeta\in\mathbb{T}.\]
The infinite sum\[
(D\widetilde{\mathfrak{H}})^{-}\oplus(D\widetilde{\mathfrak{H}})^{-}\oplus(D\widetilde{\mathfrak{H}})^{-}\oplus\cdots\]
appearing in Proposition \ref{pro:H+D+D+D_model} can then be
identified with $H^{2}((L_{\Delta}L^{2}(\mathfrak{E})^{-}))$. The
elements in this space can be viewed as functions of two variables
$(w,\zeta)\in\mathbb{D}\times\mathbb{T}$, analytic in $w$ and
measurable in $\zeta$.

We are now ready to reformulate Proposition \ref{pro:H+D+D+D_model}.
\begin{prop}
\label{pro:functional_model}Let
$\Theta:\mathbb{D}\to\mathcal{L}(\mathfrak{E})$ be a contractive
analytic function, and set
$\Delta(\zeta)=(I-\Theta(\zeta)^{*}\Theta(\zeta))^{1/2},\zeta\in\mathbb{T}.$
\begin{enumerate}
\item [{\rm(1)}]The space $(L_{\Delta}L^{2}(\mathfrak{E}))^{-}$ is reducing
for $U_{\mathfrak{E}}$.
\item [{\rm(2)}]Define the Hilbert space \[
\mathfrak{K}=H^{2}(\mathfrak{E})\oplus
H^{2}((L_{\Delta}L^{2}(\mathfrak{E}))^{-}),\] and the operators
$W_{0},W_{1}\in\mathcal{L}(\mathfrak{H})$ by \[ W_{0}(f\oplus
g)=a\oplus b,\quad W_{1}(f\oplus g)=c\oplus d,\] where\[
a(z)=zf(z),\quad b(w,\zeta)=\zeta g(w,\zeta),\]
\[
c(z)=\Theta(z)f(z),\quad d(w,\zeta)=\Delta(\zeta)f(\zeta)+wg(w,\zeta)\]
for $z,w\in\mathbb{D}$ and $\zeta\in\mathbb{T}$. Then $(W_{0},W_{1})$
is a $\{0\}$-pure bi-isometry whose characteristic function is unitarily
equivalent to $\Theta$.
\end{enumerate}
\end{prop}
We will use the notation $\mathbb{W}(\Theta)=(W_{0},W_{1})$ for the
bi-isometry described in the preceding statement. The mapping
$\Theta\mapsto\mathbb{W}(\Theta)$ establishes a bijection between
unitary equivalence classes of contractive analytic functions
$\Theta:\mathbb{D}\to\mathcal{L}(\mathfrak{E})$ and unitary
equivalence classes of $\{0\}$-pure bi-isometries $\mathbb{W}$. The
formulas given for $\mathbb{W}(\Theta)$ allow, in principle,
explicit calculations. A first instance is the following result.
\begin{prop}
\label{pro:symmetric_case}Let
$\Theta:\mathbb{D}\to\mathcal{L}(\mathfrak{E})$ be a contractive
analytic function, and denote $(W_{0},W_{1})=\mathbb{W}(\Theta)$.
\begin{enumerate}
\item [{\rm(1)}]The operator $W_{1}$ is unitary if and only if $\Theta$
is a constant unitary operator, that is, $\Theta(z)\equiv\Theta(0),$
and $\Theta(0)$ is a unitary operator in
$\mathcal{L}(\mathfrak{E})$.
\item [{\rm(2)}]The following conditions are equivalent:

\begin{enumerate}
\item $\mathbb{W}(\Theta)$ is $\{1\}$-pure.
\item the contraction $\Theta(0)$ is completely nonunitary.
\end{enumerate}
\end{enumerate}
\end{prop}
\begin{proof}
If $W_{1}$ is unitary, $W_{0}$ must be a unilateral shift, and
therefore $V=S_{\mathfrak{E}}$, and $W=T_{\Theta}$. It is well-known
that $T_{\Theta}$ is unitary if and only if $\Theta$ is a constant
unitary operator.

To prove (2), assume first that $W_{1}|\mathfrak{N}$ is unitary for
some nonzero reducing subspace $\mathfrak{N}$ of
$\mathbb{W}(\Theta)$. Applying part (1) of Proposition
\ref{cor:char_func} and part (1) of this proposition, which has
already been proved, shows that we can write
$\Theta=\Theta'\oplus\Theta''$, with $\Theta''$ a constant unitary
operator acting on a nonzero space. In particular $\Theta(0)$ has a
nontrivial unitary direct summand. Conversely, assume that
$\Theta(0)$ is not completely nonunitary, so that its restriction to
some nonzero invariant subspace$\mathfrak{E}_{0}$ is a unitary
operator. The contractive analytic function
$\Theta_{0}:\mathbb{D}\to\mathcal{L}(\mathfrak{E}_{0})$ defined by
$\Theta_{0}(z)=P_{\mathfrak{E}_{0}}\Theta(z)|\mathfrak{E}_{0}$ is
such that $\Theta_{0}(0)$ is unitary. The maximum principle implies
that $\Theta_{0}$ is constant, and $\mathfrak{E}_{0}$ reduces each
$\Theta(z)$ to $\Theta_{0}$. A second application of part (1) of
Proposition \ref{cor:char_func}, as well as the already proved part
(1) of this proposition, shows that $W_{1}|\mathfrak{N}$ is unitary
for some nonzero reducing subspace $\mathfrak{N}$ of
$\mathbb{W}(\Theta)$.
\end{proof}
If $\mathbb{W}=\mathbb{W}(\Theta)$, one can also calculate the characteristic
function of $\check{\mathbb{W}}=(W_{1},W_{0})$, whose coefficients
are\[
(I-W_{1}W_{1}^{*})W_{1}^{*k}W_{0}|\text{ran}(I-W_{1}W_{1}^{*}),\quad k\ge0.\]
Thus this function is given by\[
(I-W_{1}W_{1}^{*})(I-zW_{1}^{*})^{-1}W_{0}|\text{ran}(I-W_{1}W_{1}^{*}),\quad z\in\mathbb{D}.\]
In these formulas we use the abbreviation `ran' for the range of an
operator.

\section{The structure of bi-shifts\label{sec:bi-shifts}}

Consider a bi-isometry $\mathbb{W}=(W_{0},W_{1})$. As seen earlier,
the operators $W_{0}$ and $W_{1}$ do not need to be cnu, even if
$\mathbb{W}$ is $\{0\}$-pure and $\{1\}$-pure. In this section
we study bi-isometries for which both $W_{0}$ and $W_{1}$ are cnu,
and such bi-isometries will be called \emph{bi-shifts}. Clearly bi-shifts
are both $\{0\}$-pure and $\{1\}$-pure. Note that the bi-shifts
described in \cite{Gasp} are, in our terminology, doubly commuting
bi-shifts; see Proposition \ref{pro:D4} below.
\begin{prop}
\label{pro:D1}Assume that the bi-isometry $\mathbb{W}$ is both $\{0\}$-pure
and $\{1)$-pure. The following conditions are equivalent.
\begin{enumerate}
\item [(1)]$\mathbb{W}$ is a bi-shift.
\item [(2)]$W_{0}^{*n}\to0$ and $W_{1}^{*n}\to0$ as $n\to\infty$ in
the strong operator topology.
\item [(3)]The characteristic function $\Theta_{\mathbb{W}}$ is inner
\emph{(}that is, $\Theta_{\mathbb{W}}(\zeta)\in\mathcal{L}(\mathfrak{E})$
is an isometry for almost every $\zeta\in\mathbb{T}$\emph{) }and
it enjoys the following property:

\begin{enumerate}
\item [($*$)]There exists no inner function $\Omega:\mathbb{D}\to\mathcal{L}(\mathfrak{F},\mathfrak{E})$
such that $\mathfrak{F}\ne\{0\}$ and\[
\Theta_{\mathbb{W}}(z)\Omega(z)=\Omega(z)U,\quad z\in\mathbb{D},\]
with a unitary operator $U\in\mathcal{L}(\mathfrak{F})$.
\end{enumerate}
\end{enumerate}
\end{prop}
\begin{proof}
The proposition is almost immediate, but we provide the brief argument
below in order to illustrate the use of the results in the preceding
section. The equivalence between (1) and (2) follows from the fact
that an isometry is cnu if and only if it is a unilateral shift. Assume
next that $(2)$ holds so that, in particular, $W_{0}$ has no unitary
part. With the notation of the preceding sections, $V_{0}=W_{0}$,
and $A=W_{1}$, so that $\mathbb{W}$ serves as its own characteristic
pair. Passing to the functional model, we identify $W_{0}$ with the
unilateral shift $S_{\mathfrak{E}}$, in which case $W_{1}=T_{\Theta}$
for some operator valued function $\Theta$. The function $\Theta$
must then be inner because $T_{\Theta}$ is an isometry. Assume now
that a function $\Omega$ exists with the properties in $(*)$. Then
it follows that $T_{\Theta}|\Omega H^{2}(\mathfrak{F})$ is a unitary
operator, unitarily equivalent to $T_{U}\in\mathcal{L}(H^{2}(\mathfrak{F}))$.
This contradicts the assumption that (2) holds, and we conclude that
(3) is true. Finally, assume that (3) holds, but (2) does not. Since
$S_{\mathfrak{E}}$ is completely nonunitary, the operator $T_{\Theta}$
must have a unitary part. The nonzero space \[
\mathfrak{M}=\bigcap_{n=0}^{\infty}W_{1}^{n}\mathfrak{H}=\bigcap_{n=0}^{\infty}T_{\Theta}^{n}H^{2}(\mathfrak{E})\]
on which this unitary part acts is obviously invariant for $S_{\mathfrak{E}}$,
and the Beurling-Lax-Halmos theorem implies that $\mathfrak{M}=\Omega H^{2}(\mathfrak{F})$
for some inner function $\Omega:\mathbb{D}\to\mathcal{L}(\mathfrak{F},\mathfrak{E})$
with $\mathfrak{F}\ne\{0\}$. The operator $T_{\Omega}^{-1}T_{\Theta}T_{\Omega}$
is then a unitary operator in the commutant of $S_{\mathfrak{F}}$,
and such operators are of the form $T_{U}$ for some unitary operator
$U\in\mathcal{L}(\mathfrak{F})$. We conclude that $T_{\Theta}T_{\Omega}=T_{\Omega}T_{U}$,
contrary to (3).\end{proof}
\begin{rem}
The example $\Theta(z)\equiv I_{\mathfrak{E}}$, $z\in\mathbb{D}$,
shows that condition $(*)$ is needed in the third statement of the
preceding proposition.
\end{rem}
The preceding result shows that the construction of bi-shifts requires
the construction of appropriate inner functions $\Theta:\mathbb{D}\to\mathcal{L}(\mathfrak{E})$.
We start with some simple examples. Fix a nonzero Hilbert space $\mathfrak{E}$
and an inner function $\vartheta\in H^{\infty}$. We can then form
the bi-isometry $\mathbb{W}(\vartheta\otimes I_{\mathfrak{E}})=(S_{\mathfrak{E}},T_{\vartheta\otimes I_{\mathfrak{E}}})$.
This is easily seen to be a bi-shift provided that $\vartheta$ is
not constant.
\begin{prop}
\label{pro:D3}Given two nonconstant inner functions $\vartheta_{1},\vartheta_{2}\in H^{\infty}$,
the bi-shifts $\mathbb{W}(\vartheta_{1}\otimes I_{\mathfrak{E}})$
and $\mathbb{W}(\vartheta_{2}\otimes I_{\mathfrak{E}})$ are quasi-similar
if and only if $\vartheta_{1}=\vartheta_{2}$.\end{prop}
\begin{proof}
Let $X\in\mathcal{I}(\mathbb{W}(\vartheta_{1}\otimes I_{\mathfrak{E}}),\mathbb{W}(\vartheta_{2}\otimes I_{\mathfrak{E}}))$
be a quasi-affinity. We have $X\in(S_{\mathfrak{E}})'$, and therefore
we have $X=T_{\Xi}$ for some outer function $\Xi:\mathbb{D}\to\mathcal{L}(\mathfrak{E})$.
The relation $XT_{\vartheta_{1}\otimes I_{\mathfrak{E}}}=T_{\vartheta_{2}\otimes I_{\mathfrak{E}}}X$
implies then\[
(\vartheta_{1}(z)-\vartheta_{2}(z))\Xi(z)=\Xi(z)(\vartheta_{1}(z)\otimes I_{\mathfrak{E}})-(\vartheta_{2}(z)\otimes I_{\mathfrak{E}})\Xi(z),\quad z\in\mathbb{D}.\]
The operator $\Xi(z)$ has dense range for every $z\in\mathbb{D}$,
and we conclude that $\vartheta_{1}=\vartheta_{2}$. The converse
is immediate.
\end{proof}
As pointed out earlier, the bi-isometries $\mathbb{W}(\Theta_{1})$
and $\mathbb{W}(\Theta_{2})$, $\Theta_{1},\Theta_{2}:\mathbb{D}\to\mathfrak{E}$,
are unitarily equivalent if and only if the functions $\Theta_{1}$
and $\Theta_{2}$ are unitarily equivalent, that is, $U\Theta_{1}(z)=\Theta_{2}(z)U$
for a unitary operator $U$ independent of $z\in\mathbb{D}$. Similarity
of the two bi-isometries requires the existence of an invertible outer
function $\Psi:\mathbb{D}\to\mathcal{L}(\mathfrak{E})$ such that
$\Psi(z)\Theta_{1}(z)=\Theta_{1}(z)\Psi(z)$ for all $z\in\mathbb{D}$.

Another important family of bi-shifts is defined on the Hardy space
$H^{2}(\mathbb{D}^{2})\otimes\mathfrak{F}$ by the formula $\mathbb{W}_{\mathfrak{F}}=(W_{0},W_{1})$,
where \[
(W_{j}f)(z_{0},z_{1})=z_{j}f(z_{0},z_{1}),\quad f\in H^{2}(\mathbb{D}^{2})\otimes\mathfrak{F},(z_{0},z_{1})\in\mathbb{D}^{2}.\]
This class of bi-isometries has a simple characterization. Parts of
the following proposition are known. We include a brief argument for
the reader's convenience.
\begin{prop}
\label{pro:D4}Assume that the bi-isometry $\mathbb{W}$ is both $\{0\}$-pure
and $\{1\}$-pure. The following conditions are equivalent.
\begin{enumerate}
\item [(1)]$\mathbb{W}$ is unitarily equivalent to $\mathbb{W}_{\mathfrak{F}}$
for some Hilbert space $\mathfrak{F}$.
\item [(2)]$\mathbb{W}$ is doubly commuting, that is, $W_{0}W_{1}^{*}=W_{1}^{*}W_{0}$.
\item [(3)]The characteristic function $\Theta_{\mathbb{W}}$ is a constant
isometry.
\item [(4)]The pivotal operator of $(W_{1},W_{0})$ is an isometry.
\item [(5)] The pivotal operator of $\mathbb{W}$ is an isometry.
\end{enumerate}
\end{prop}
\begin{proof}
It is immediate that $(1)$ implies (2). For the remainder of the
argument we identify $\mathbb{W}$ with $\mathbb{W}(\Theta),$ where
$\Theta:\mathbb{D}\to\mathcal{L}(\mathfrak{E})$ is a contractive
analytic function. Thus $\mathbb{W}$ acts on the space $\mathfrak{K}$
described in Proposition \ref{pro:symmetric_case}. Assume now that
(2) holds. In this case the kernel of $W_{0}^{*}$ must be a reducing
subspace for $W_{1}$. This kernel consists of functions in $\mathfrak{K}$
of the form $e\oplus0\oplus0\oplus\cdots$, with $e\in\mathfrak{E}$
a constant. Since \[
W_{1}(e\oplus0\oplus0\oplus\cdots)=\Theta e\oplus\Delta e\oplus0\oplus\cdots,\]
 we deduce immediately that $\Theta$ is constant and $\Delta=0$,
so that $(3)$ is true. Assume now that (3) holds, so that $\Theta$
is a constant isometry. It follows that $\Theta(0)$ is in particular
an isometry. Condition (4) follows because $\Theta(0)$ is the pivotal
operator of the pair $(W_{1},W_{0})$. Assume that (4) holds, so that
$\Theta(0)$ is an isometry. Then it follows from the maximum principle
that $\Theta(z)=\Theta(0)$ for all $z$. In particular, the function
$\Theta$ is inner, and hence $W_{0}=S_{\mathfrak{E}}$ and $W_{1}=T_{\Theta}$.
Note that any orthogonal decomposition $\Theta(0)=\Theta_{1}\oplus\Theta_{2}$
yields a decomposition $T_{\Theta}=T_{\Theta_{1}}\oplus T_{\Theta_{2}}$.
If $\Theta_{1}$ is unitary, the operator $T_{\Theta_{1}}$ is unitary
as well, and therefore $\Theta_{1}$ must act on the space $\{0\}$
because $\mathbb{W}$ was assumed to be $\{1\}$-pure. We deduce that
$\Theta(0)$ is cnu, and thus it is unitarily equivalent to $S_{\mathfrak{F}}$
for some Hilbert space $\mathfrak{F}$, and in this case $\mathbb{W}(\Theta)$
is unitarily equivalent to $\mathbb{W}_{\mathfrak{F}}$.

So far we have proved that conditions (1--4) are equivalent. The equivalence
of (5) with these conditions follows from the symmetry of (2).
\end{proof}
The example of the constant function $\Theta(z)\equiv I$, $z\in\mathbb{D}$,
shows why the assumption that $\mathbb{W}$ is both $\{0\}$-pure
and $\{1\}$-pure is needed in the preceding proposition.

If two isometries are quasi-similar and one of them is a shift, then
the other one is a shift as well. It follows that a bi-isometry quasi-similar
to a bi-shift must also be a bi-shift. We conclude this section with
some simple properties of those bi-shifts which are similar to $\mathbb{W}_{\mathfrak{F}}$
for some $\mathfrak{F}$.
\begin{prop}
\label{pro:D5-} Let $\mathbb{W}=\mathbb{W}(\Theta)$ be a bi-shift,
where $\Theta:\mathbb{D}\to\mathcal{L}(\mathfrak{E})$ is an inner
analytic function. Assume further that $\mathbb{W}$ is similar to
$\mathbb{W}_{\mathfrak{F}}$ for some Hilbert space $\mathfrak{F}$.
Then the following assertions are true.
\begin{enumerate}
\item [(1)]The pivotal operator is similar to a unilateral shift.
\item [(2)]There exists a bounded analytic function $\Omega:\mathbb{D}\to\mathcal{L}(\mathfrak{E})$
such that \[ \Omega(z)\Theta(z)=I,\quad z\in\mathbb{D}.\]

\item [(3)]The operator $\Theta(z)$ is similar to a unilateral shift for
every $z\in\mathbb{D}$.\end{enumerate}
\begin{proof}
We argue first that two similar bi-isometries have similar pivotal
operators. Indeed, assume that
$X\in\mathcal{I}(\mathbb{W}^{(1)},\mathbb{W}^{(2)})$ is an
invertible operator. We have then $X\ker W_{0}^{(1)*}=\ker
W_{0}^{(2)*}$, and this implies that $X|\ker W_{0}^{(1)*}$ is an
invertible operator intertwining the two pivotal operators. Now, the
pivotal operator of $\mathbb{W}_{\mathfrak{F}}$ is a shift, and the
preceding observation implies (1). By symmetry, we also deduce that
$\Theta(0)$ is similar to a shift, and then (2) follows from the
main result of \cite{NF-sim}. To verify (3), we observe that the
bi-shift $\mathbb{W}_{\mathfrak{F}}$ is unitarily equivalent to
$\mathbb{W}(\Theta_{1})$, where $\Theta_{1}(z)\equiv S$ for
$z\in\mathbb{D}$, with $S\in\mathcal{L}(\mathfrak{E})$ a unilateral
shift. Let
$X\in\mathcal{I}(\mathbb{W}(\Theta),\mathbb{W}(\Theta_{1}))$ be an
invertible operator. We have $X\in(W_{0})'$, and therefore the
operator $X$ is of the form $X=T_{\Xi}$ for some bounded analytic
function $\Xi\in\mathcal{L}(\mathfrak{E})$. The fact that $X$ is
invertible implies that $X(z)$ is invertible for every
$z\in\mathbb{D}$, and the relation $XT_{\Theta}=T_{\Theta_{1}}X$
shows that $\Theta(z)$ is similar to $S=\Theta_{1}(z)$. The
proposition is proved.
\end{proof}
\end{prop}
A different approach to the similarity between a contraction and an
isometry is described in \cite{kwon}. This approach may also be
useful in the study of similarities between bi-shifts.

Conditions (1) and (2) in the above proposition are not sufficient
to imply the similarity of $\mathbb{W}(\Theta)$ to a bi-shift of the
form $\mathbb{W}_{\mathfrak{F}}$, as shown by the following example.
\begin{example}
Define $\Theta(z)\in\mathcal{L}(\ell^{2})$ using the infinite
matrix\[ \Theta(z)=\left[\begin{array}{ccccc}
\frac{3}{5}\varphi(z) & 0 & 0 & 0 & \cdots\\
\frac{4}{5}z & 0 & 0 & 0 & \cdots\\
0 & 1 & 0 & 0 & \cdots\\
0 & 0 & 1 & 0 & \cdots\\
\vdots & \vdots & \vdots & \ddots & \cdots\end{array}\right],\quad
z\in\mathbb{D},\] where $\varphi\in H^{\infty}$ is an inner function
such that $\varphi(0)\ne0$. The operator $\Theta(0)$ has the
eigenvalue $3\varphi(0)/5$ and therefore it is not similar to a
shift. However $\Theta$ satisfies condition (2) in the preceding
proposition. One left inverse is given by\[
\Omega(z)=\left[\begin{array}{ccccc}
\frac{5}{3\varphi(0)} & \eta(z) & 0 & 0 & \cdots\\
0 & 0 & 1 & 0 & \cdots\\
0 & 0 & 0 & 1 & \cdots\\
0 & 0 & 0 & 0 & \ddots\\
\vdots & \vdots & \vdots & \vdots & \vdots\end{array}\right],\quad
z\in\mathbb{D},\] with \[
\eta(z)=\frac{4}{5z}\left[1-\frac{\varphi(z)}{\varphi(0)}\right],\quad
z\in\mathbb{D}.\] The reader will verify without difficulty that
$\mathbb{W}(\Theta)$ is indeed a bi-shift.
\end{example}

\section{The unitary invariants of a functional model}

Bi-isometries $\mathbb{W}=(W_{0},W_{1})$ with the property that the
product $W_{0}W_{1}$ is a shift were classified, up to unitary
equivalence, in \cite{BCL}; see also \cite{BDF}. The parameters in
that classification are pairs $(U,P)$, where $U$ is a unitary
operator on a Hilbert space $\mathfrak{D}$, and $P$ is an orthogonal
projection on $\mathfrak{D}$. In this section we consider the
characteristic functions of such bi-isometries. The bi-isometry
$\mathbb{W}=(W_{0},W_{1})$ associated to the pair $(U,P)$ acts on
$H^{2}(\mathfrak{D})$ and is defined by\[
(W_{0}f)(z)=U(zP+P^{\perp})f(z),\enskip(W_{1}f)(z)=(P+zP^{\perp})U^{*}f(z),\quad
f\in H^{2}(\mathfrak{D}),z\in\mathbb{D}.\] The space $\mathfrak{D}$
is identified with the space $\ker(W_{0}W_{1})^{*}$ of constant
functions in $H^{2}(\mathfrak{D})$, while the range of $P^{\perp}$
is identified with $\ker W_{1}^{*}$. For a constant function
$f_{0}\in\mathfrak{D}$ we have \begin{equation}
Uf_{0}=W_{0}f_{0},\quad f_{0}\in
P^{\perp}\mathfrak{D},\label{eq:Uf0-first}\end{equation} while for
$f_{0}\in P\mathfrak{D}$ we have \[
W_{0}f_{0}=zUf_{0}=W_{0}W_{1}Uf_{0}.\] Therefore the vector
$f_{0}=W_{1}Uf_{0}$ is in the range of $W_{1}$, and we find that
\begin{equation} Uf_{0}=W_{1}^{*}f_{0},\quad f_{0}\in
P\mathfrak{D}.\label{eq:Uf_0-second}\end{equation} From this we
easily conclude that\[ \ker
W_{0}^{*}=UP\mathfrak{D}=W_{1}^{*}P\mathfrak{D}.\] By reversing the
order of these observations we easily deduce the following result.
\begin{prop}
\label{pro:D=00003DWE+F}Let $\mathbb{W}=(W_{0},W_{1})$ be a
bi-isometry on $\mathfrak{H}$. Define spaces\[
\mathfrak{D}=\ker(W_{0}W_{1})^{*},\quad\mathfrak{E}=\ker
W_{0}^{*},\quad\mathfrak{F}=\ker W_{1}^{*}.\]

\begin{enumerate}
\item [{\rm(1)}]We have $\mathfrak{D}=\mathfrak{E}\oplus W_{0}\mathfrak{F}=W_{1}\mathfrak{E}\oplus\mathfrak{F}$.
\item [{\rm(2)}]The operator $U:\mathfrak{D}\to\mathfrak{D}$ defined by\[
U(W_{1}e+f)=e+W_{0}f\quad e\in\mathfrak{E},f\in\mathfrak{F},\] is
unitary.
\item [{\rm(3)}]The bi-isometry associated with the pair $(U,P_{W_{1}\mathfrak{E}})$
on $\mathfrak{D}$ is unitarily equivalent to the cnu part of
$\mathbb{W}$.
\end{enumerate}
\end{prop}
For further calculation, it is convenient to replace the space
$\mathfrak{D}$ by the external direct sum
$\mathfrak{E}\oplus\mathfrak{F}$ via the identification
$\Phi:e\oplus f\mapsto W_{1}e+f$. With this identification we
obviously have\[ \Phi^{*}P\Phi=\left[\begin{array}{cc}
I_{\mathfrak{E}} & 0\\
0 & 0\end{array}\right].\]

\begin{cor}
\label{cor:(U,P)on(E+F)}With the notation of Proposition \emph{\ref{pro:D=00003DWE+F}},
we have\[
\Phi^{*}U\Phi=\left[\begin{array}{cc}
W_{1}^{*}|\mathfrak{E} & W_{1}^{*}W_{0}|\mathfrak{F}\\
(I-W_{1}W_{1}^{*})|\mathfrak{E} &
(I-W_{1}W_{1}^{*})W_{0}|\mathfrak{F}\end{array}\right].\]
\end{cor}
\begin{proof}
For a vector $e\in\mathfrak{E}$ we have\[
U\Phi(e\oplus0)=UW_{1}e=e=W_{1}W_{1}^{*}e+(I-W_{1}W_{1}^{*})e,\] and
this is precisely the decomposition of this vector as an element of
the space $W_{1}\mathfrak{E}\oplus\mathfrak{F}$. Therefore\[
\Phi^{*}U\Phi(e\oplus0)=W_{1}^{*}e\oplus(I-W_{1}W_{1}^{*})e.\] To
verify the identity involving the second column, we use a similar
calculation:\[ U\Phi(0\oplus
f)=Uf=W_{0}f=W_{1}W_{1}^{*}W_{0}f+(I-W_{1}W_{1}^{*})W_{0}f,\quad
f\in\mathfrak{F}.\]
 In these calculations we made use of (\ref{eq:Uf0-first}) and (\ref{eq:Uf_0-second}).
\end{proof}
Let us consider now a contractive analytic function
$\Theta:\mathbb{D}\to\mathcal{L}(\mathfrak{E})$ and the functional
model $\mathbb{W}(\Theta)=(W_{0},W_{1})$. In order to identify the
space $\mathfrak{F}$, it will be useful to recall a few facts from
the theory of functional models of contraction operators. Let us
introduce the auxiliary space\[
\mathfrak{K}=H^{2}(\mathfrak{E})\oplus(L_{\Delta}L^{2}(\mathfrak{E}))^{-},\]
which can be viewed as a subspace of
$\mathfrak{H}=H^{2}(\mathfrak{E})\oplus
H^{2}((L_{\Delta}L^{2}(\mathfrak{E}))^{-})$. Obviously, the space
$\mathfrak{K}$ is reducing for $W_{0}$. The space\[
\mathfrak{G}=\{\Theta u\oplus\Delta u:u\in H^{2}(\mathfrak{E})\}\]
is invariant for $W_{0}$, and therefore \[
\mathfrak{H}(\Theta)=\mathfrak{K}\ominus\mathfrak{G}\] is invariant
for $W_{0}^{*}$. The compression of $W_{0}$ to this space is denoted
$S(\Theta)$, and it is called the functional model associated with
$\Theta$. It is known that $S(\Theta)$ is a completely nonunitary
contraction, and the characteristic function of $S(\Theta)$
coincides (in the sense defined in \cite{SzNFbook-1}) with the
purely contractive part of the function $\Theta$.

A vector $u\oplus v\in\mathfrak{K}$ belongs to
$\mathfrak{H}(\Theta)$ if and only if the measurable function
$\Theta^{*}u+\Delta v$ is orthogonal to $H^{2}(\mathfrak{E})$. In
other words, we have a Fourier expansion\[ \Theta^{*}u+\Delta
v=\sum_{n=-1}^{\infty}\zeta^{n}e_{n},\] with $e_{n}\in\mathfrak{E}$.
We will use the notation $(\Theta^{*}u+\Delta v)_{-1}$ for $e_{-1}$.
\begin{lem}
\label{lem:H(teta)=00003DF}Viewed as a subspace of $\mathfrak{H}$,
we have $\mathfrak{H}(\Theta)=\mathfrak{F}$. Moreover, $S(\Theta)$
is precisely the pivotal operator associated with the bi-isometry
$\mathbb{W}(\Theta)$:\[ S(\Theta)^{*}=W_{0}^{*}|\mathfrak{F}.\]
\end{lem}
\begin{proof}
In order to identify $\mathfrak{F}$, we consider its orthogonal
complement which is easily calculated as\[
\mathfrak{F}^{\perp}=W_{1}\mathfrak{H}=\mathfrak{G}\oplus
W_{1}H^{2}((L_{\Delta}L^{2}(\mathfrak{E}))^{-}).\] The conclusion
$\mathfrak{H}(\Theta)=\mathfrak{F}$ then follows because\[
\mathfrak{H}=\mathfrak{K}\oplus
W_{1}H^{2}((L_{\Delta}L^{2}(\mathfrak{E}))^{-}).\] The
identification of the pivotal operator follows now from the fact
that $\mathfrak{H}(\Theta)=\mathfrak{F}$ is invariant for
$W_{0}^{*}$.\end{proof}
\begin{prop}
\label{pro:uni_invariants_V(teta)}Let
$\Theta:\mathbb{D}\to\mathcal{L}(\mathfrak{E})$ be a contractive
analytic function, and $\mathbb{W}(\Theta)=(W_{0},W_{1})$ the
corresponding model bi-isometry. Then $\mathbb{W}(\Theta)$ is
unitarily equivalent to the bi-isometry associated with the pair
$(U,P)$ of operators on $\mathfrak{E}\oplus\mathfrak{H}(\Theta)$
defined as follows:\begin{eqnarray}
U(e\oplus0) & = & \Theta(0)^{*}e\oplus[(e-\Theta\Theta(0)^{*}e)\oplus(-\Delta\Theta(0)^{*}e)],\quad e\in\mathfrak{E,}\nonumber \\
U(0\oplus(u\oplus v)) & = & (\Theta^{*}u+\Delta v)_{-1}\oplus
S(\Theta)(u\oplus v),\quad u\oplus
v\in\mathfrak{H}(\Theta),\label{eq:prop5-4}\end{eqnarray} and\[
P=\left[\begin{array}{cc}
I_{\mathfrak{E}} & 0\\
0 & 0\end{array}\right].\]
\end{prop}
\begin{proof}
This proof amounts to an identification of the matrix entries in
Corollary \ref{cor:(U,P)on(E+F)}. It is convenient to regard
$\mathfrak{H}$ as an infinite orthogonal sum\[
\mathfrak{H}=H^{2}(\mathfrak{E})\oplus(L_{\Delta}L^{2}(\mathfrak{E}))^{-}\oplus(L_{\Delta}L^{2}(\mathfrak{E}))^{-}\oplus\cdots,\]
relative to which the operator $W_{1}$ has the matrix\[
W_{1}=\left[\begin{array}{cccc}
T_{\Theta} & 0 & 0 & \cdots\\
L_{\Delta}|H^{2}(\mathfrak{E}) & 0 & 0 & \cdots\\
0 & I_{(L_{\Delta}L^{2}(\mathfrak{E}))^{-}} & 0 & \cdots\\
\vdots & \vdots & \vdots & \ddots\end{array}\right].\] We now apply
the formulas in Corollary \ref{cor:(U,P)on(E+F)} to calculate the
entries of the matrix $U$ explicitly. Thus, for $e\in\mathfrak{E}$,
which is viewed now as a subspace of $\mathfrak{H}$, we obtain by
applying the matrix above\[
W_{1}^{*}e=T_{\Theta}^{*}e=P_{H^{2}(\mathfrak{E})}\Theta^{*}e=\Theta(0)^{*}e,\]
and \[ (I-W_{1}W_{1}^{*})e=e-W_{1}\Theta(0)^{*}e.\]
 If $u\oplus v\in\mathfrak{H}(\Theta)$ then clearly\[
(I-W_{1}W_{1}^{*})W_{0}(u\oplus
v)=P_{\mathfrak{H}(\Theta)}W_{0}(u\oplus v)=S(\Theta)(u\oplus v).\]
For the first direct summand in the right-hand side of
(\ref{eq:prop5-4}), let us write $W_{0}(u\oplus v)=u'\oplus v'$ and
note that\[ W_{1}^{*}W_{0}(u\oplus v)=W_{1}^{*}(u'\oplus
v')=P_{H^{2}(\mathfrak{E})}(\Theta^{*}u'+\Delta v').\] If we write
the Fourier expansion\[ \Theta^{*}u+\Delta
v=\sum_{n=-1}^{\infty}\zeta^{n}e_{n},\] then\[ \Theta^{*}u'+\Delta
v'=\sum_{n=-1}^{\infty}\zeta^{n+1}e_{n},\] and the projection of
this function onto $H^{2}(\mathfrak{E})$ is precisely
$e_{-1}=(\Theta^{*}u+\Delta v)_{-1}$, as stated.
\end{proof}

\section{Examples of irreducible two-isometries and direct integral decompositions\label{sec:Examples-of-irreducible}}

For a single isometry, that is, when $I$ has only one element, it
follows from the von Neumann--Wold theorem that there is, up to
unitary equivalence, only one nonunitary irreducible isometry.
However, when $I$ has two or more elements there are many
irreducible families of commuting isometries which do not consist of
unitary operators. We will illustrate this in the case of
bi-isometries $\mathbb{V}=(V_{0},V_{1})$. We recall that a complete
unitary invariant of a completely nonunitary bi-isometry
$\mathbb{V}=(V_{0},V_{1})$ is given by a pair $(U,P)$, where $U$ is
a unitary operator on some Hilbert space $\mathfrak{D}$, and $P$ is
an orthogonal projection on $\mathfrak{D}$. The bi-isometry
determined by $(U,P)$ acts on $H^{2}(\mathfrak{D})$ as follows:\[
(V_{0}f)(z)=U(zP+P^{\perp})f(z),\quad(V_{1}f)(z)=(P+zP^{\perp})U^{*}f(z),\quad
f\in H^{2}(\mathfrak{D}),z\in\mathbb{D},\] where
$P^{\perp}=I_{\mathfrak{E}}-P$. The bi-isometry $\mathbb{V}$ is
irreducible if and only if the pair $(U,P)$ is irreducible. Note for
further use that the product $V_{0}V_{1}$ is precisely
multiplication by the variable $z$. (These unitary invariants
classify more general bi-isometries than the completely nonunitary
ones; see \cite{BCL,BDF}.)

For our illustration we will let $U$ be the bilateral shift on the
space $L^{2}$ of all square integrable functions on the unit circle
$\mathbb{T}$; thus\[
(Uf)(\zeta)=\zeta f(\zeta),\quad f\in L^{2},\zeta\in\mathbb{T}.\]
 We will denote by $e_{j}(\zeta)=\zeta^{j}$ the standard orthonormal
basis in $L^{2}$, and for every set $A\subset\mathbb{Z}$ of integers
we denote by $Q_{A}$ the orthogonal projection onto the space generated
by $\{e_{j}:j\in A\}$. In this case $V_{0}$ and $V_{1}$ are uniquely
determined by the relations $V_{1}e_{n+1}=e_{n}$ if $n\in A$ and
$V_{0}e_{n}=e_{n+1}$ if $n\notin A$.
\begin{prop}
\label{pro:equivalence=00003Dtranslation}Two pairs $(U,Q_{A})$,
$(U,Q_{B})$ are unitarily equivalent if and only if there exists
$n\in\mathbb{Z}$ such that \[
B=\{i+n:i\in A\}.\]
\end{prop}
\begin{proof}
Sufficiency is obvious: if $B=A+n$ then the operator $U^{n}$ implements
the unitary equivalence of the two pairs. Conversely, assume that
there is a unitary operator $\Phi$ on $L^{2}$ such that $\Phi U=U\Phi$
and $UQ_{A}=Q_{B}U$. There exists then a function $\varphi\in L^{\infty}$
such that $|\varphi|=1$ almost everywhere and $\Phi f=\varphi f$
for every $f\in L^{2}$. The fact that $\varphi e_{i}$ is in the
range of $Q_{B}$ for $i\in A$ means that\[
(\varphi,e_{j-i})=(\varphi e_{i},e_{j})=0,\quad i\in A,j\notin B.\]
Similarly, $\varphi e_{i}$ is in the range of $Q_{B}^{\perp}$ if
$i\notin A$, so that\[
(\varphi,e_{j-i})=0,\quad i\notin A,j\in B.\]
We deduce that there exists at least one integer $n$ not in the set
$\{j-i:(i,j)\in(A\times(\mathbb{Z}\setminus B))\cup((\mathbb{Z}\setminus A)\times B\}$.
The function $e_{n}$ will then have the property that $e_{n+i}=e_{n}e_{i}$
is in the range of $Q_{B}$ if $i\in A$, and it is in the range of
$Q_{B}^{\perp}$ if $i\notin A$. Therefore $B=A+n$.\end{proof}
\begin{cor}
\label{cor:irreducible(U,Q)}The pair $(U,Q_{A})$ is reducible if
and only of $A$ is a periodic set, that is, $A=A+n$ for some nonzero
integer $n$.\end{cor}
\begin{proof}
The pair $(U,Q_{A})$ is reducible if and only if it commutes with
a unitary which is not a scalar multiple of the identity. The argument
in the proof of the preceding proposition shows that such a unitary
can be chosen to be multiplication by $e_{n}$ for some $n\in\mathbb{Z}\setminus\{0\}$.
\end{proof}
We see therefore that there is a continuum of mutually inequivalent
irreducible bi-isometries. Indeed, there is a continuum of subsets
of $\mathbb{Z}$, and only countably many of them are periodic.

Quite interestingly, the bi-isometry associated with $(U,Q_{A})$ can
be described very explicitly. Consider the space
$L^{2}(\mathbb{T}^{2})=L^{2}\otimes L^{2}$ and its standard
orthonormal basis \[
e_{ij}(\zeta_{0},\zeta_{1})=\zeta_{0}^{i}\zeta_{1}^{j},\quad
i,j\in\mathbb{Z},\zeta_{0},\zeta_{1}\in\mathbb{T}.\] Multiplication
by the two variables defines a bi-isometry
$\mathbb{V}=(V_{0},V_{1})$ on $L^{2}(\mathbb{T}^{2})$; actually
$V_{0}$ and $V_{1}$ are unitary. We will look at proper nonempty
subsets $\Gamma\subset\mathbb{Z}^{2}$ with the property that the
space $\mathfrak{H}_{\Gamma}$ generated by
$\{e_{ij}:(i,j)\in\Gamma\}$ is invariant for $\mathbb{V}$. In other
words, $(i+n,j+m)\in\Gamma$ if $(i,j)\in\Gamma$ and $n,m\ge0$ or,
equivalently, $\Gamma+\mathbb{N}^{2}\subset\Gamma$. We define the
boundary $\partial\Gamma$ of $\Gamma$ to consist of those pairs
$(i,j)\in\Gamma$ such that $(i-1,j-1)$ does not belong to $\Gamma$.
For each integer $n$, there exists a unique point
$\gamma_{n}=(i_{n},j_{n})\in\partial\Gamma$ such that
$i_{n}-j_{n}=n$. Uniqueness is obvious by the definition of
$\partial\Gamma$; existence follows from the fact that
$\varnothing\ne\Gamma\ne\mathbb{Z}^{2}$. The difference
$\gamma_{n+1}-\gamma_{n}=(i_{n+1}-i_{n},j_{n+1}-j_{n})$ is either
$(1,0)$ or $(0,-1)$. We can then define the set
$A_{\Gamma}\subset\mathbb{Z}$ by\[
A_{\Gamma}=\{n\in\mathbb{Z}:\gamma_{n+1}-\gamma_{n}=(0,-1)\}.\]
 Geometrically, $A_{\Gamma}$ is the union of the vertical segments
in $\partial\Gamma$, omitting the lower endpoint of each one. The
following result is an easy exercise.
\begin{prop}
\label{pro:Gamma-vs-A}For every subset $A\subset\mathbb{Z}$ there
exists a nonempty subset $\Gamma\subset\mathbb{Z}^{2}$ such that
$\Gamma+\mathbb{N}^{2}\subset\Gamma$ and $A_{\Gamma}=A$. We have
$A_{\Gamma+(p,q)}=A_{\Gamma}+p-q$ for all $(p,q)\in\mathbb{Z}^{2}$.
\begin{prop}
\label{pro:V|HGamma-is(U,Q)}Let $\Gamma$ be a nonempty proper subset
of $\mathbb{Z}^{2}$ such that $\mathfrak{H}_{\Gamma}$ is invariant
for $\mathbb{V}$. The bi-isometry associated with the invariants
$(U,Q_{A_{\Gamma}})$ is unitarily equivalent to
$\mathbb{V}|\mathfrak{H}_{\Gamma}$.
\end{prop}
\end{prop}
\begin{proof}
The space
$\mathfrak{H}_{\partial\Gamma}=\mathfrak{H}_{\Gamma}\ominus
VW\mathfrak{H}_{\Gamma}$ can be identified with $L^{2}$ by mapping
$e_{\gamma_{n}}$ to $e_{n}$. Denote by $U_{0}$ the unitary operator
on $\mathfrak{H}_{\partial\Gamma}$ which corresponds to the shift on
$L^{2}$; in other words, $U_{0}e_{\gamma_{n}}=e_{\gamma_{n+1}}$.
Since $V_{0}V_{1}$ corresponds with multiplication by $z$, it is
clear that $\mathfrak{H}_{\partial\Gamma}$ can be identified with
$H^{2}(\mathfrak{H}_{\partial\Gamma})$. Therefore, we only need to
show that $V_{0}e_{\gamma_{n}}=V_{0}V_{1}e_{\gamma_{n+1}}$ if $n\in
A_{\Gamma}$ and $V_{0}e_{\gamma_{n}}=e_{\gamma_{n+1}}$ if $n\notin
A_{\Gamma}$. This however is immediate from the definition of
$A_{\Gamma}$ and the remark preceding Proposition
\ref{pro:equivalence=00003Dtranslation}.
\end{proof}
A direct consequence of this proposition is the following:
\begin{cor}
\label{cor:equivalence_of_restricted_bishifts}Let $\Gamma$ and
$\Gamma'$ be two nonempty proper subsets of $\mathbb{Z}^{2}$ such
that $\mathfrak{H}_{\Gamma}$ and $\mathfrak{H}_{\Gamma'}$ are
invariant for $\mathbb{V}$.
\begin{enumerate}
\item [{\rm(1)}]The bi-isometries $\mathbb{V}|\mathfrak{H}_{\Gamma}$ and $\mathbb{V}|\mathfrak{H}_{\Gamma'}$
are unitarily equivalent if and only if $\Gamma'=\Gamma+\gamma$ for
some $\gamma\in\mathbb{Z}^{2}$.
\item [{\rm(2)}]The bi-isometry $\mathbb{V}|\mathfrak{H}_{\Gamma}$ is reducible
if and only if $\Gamma=\Gamma+\gamma$ for some $\gamma\in\mathbb{Z}^{2}\setminus\{(0,0)\}$.
\end{enumerate}
\end{cor}
Two particular sets $\Gamma$ yielding irreducible bi-isometries were
considered in \cite{Gasp,DPOpov3}. The first is
$\Gamma=\mathbb{N}^{2}$, for which $A_{\Gamma}=\{n:n<0\}$. The
restriction $\mathbb{V}|\mathfrak{H}_{\Gamma}$ is a doubly commuting
bi-shift. The second is
$\Gamma=(\mathbb{Z}\times\mathbb{N})\cup(\mathbb{N}\times\mathbb{Z})$,
for which $A_{\Gamma}=\mathbb{N}$. The corresponding restriction of
$\mathbb{V}$ was called a modified bi-shift in these works. The
modified bi-shift can be seen to be the dual of the doubly commuting
 bi-shift in the sense of \cite{Conw}.

The bi-isometries of the form $\mathbb{V}|\mathfrak{H}_{\Gamma}$
were considered earlier in \cite{Salas}. They have the special property
that the range projections of the isometries in the multiplicative
semigroup they generate commute with each other. The case $\Gamma\subset\mathbb{N}^{2}$
was also considered in \cite{Dou-Nak-Seto} from the point of view
of Hilbert modules over the bidisk algebra.

We now illustrate the decomposition of a bi-isometry into a direct
integral of irreducibles with the particular case provided by the
set $A=2\mathbb{Z}$. In this case, the commutant of the pair $(U,Q_{A})$
is the algebra generated by $U^{2}$, and this operator is a unitary
operator with uniform multiplicity 2 relative to the usual arclength
measure on $\mathbb{T}$. This is realized upon using the identification
\[
\Phi:L^{2}\oplus L^{2}\to L^{2}\]
defined by\[
(\Phi(f\oplus g))(\zeta)=f(\zeta^{2})+\zeta g(\zeta^{2}),\quad\zeta\in\mathbb{T}.\]
The operator $\Phi^{*}U\Phi$ is simply multiplication by the matrix-valued
function\[
U_{0}(\zeta)=\left[\begin{array}{cc}
0 & \zeta\\
1 & 0\end{array}\right],\quad\zeta\in\mathbb{T},\]
while $\Phi^{*}Q_{A}\Phi$ is multiplication by the constant matrix\[
P_{0}=\left[\begin{array}{cc}
1 & 0\\
0 & 0\end{array}\right].\]
In other words, we have the decomposition\[
(U,P)=\int_{\mathbb{T}}^{\oplus}(U_{0}(\zeta),P_{0})|d\zeta|,\]
and it is clear that the pairs $(U_{0}(\zeta),P_{0})$ are irreducible
and mutually inequivalent. This corresponds with a direct integral
decomposition of the corresponding bi-isometry. The reader will have
no difficulty verifying that the bi-isometry associated with $(U_{0}(\zeta),P_{0})$
is of the form $(\zeta S,S)$, where $S$ is a unilateral shift of
multiplicity one.

The general case of a set $A$ such that $A=A+n$, with $n>2$, lends
itself to a similar analysis, with\[
U_{0}(\zeta)=\left[\begin{array}{cccccc}
0 & 0 & 0 & \cdots & 0 & \zeta^{n-1}\\
1 & 0 & 0 & \cdots & 0 & 0\\
0 & \zeta & 0 & \cdots & 0 & 0\\
\vdots & \vdots & \vdots & \ddots & \vdots & \vdots\\
0 & 0 & 0 & \cdots & 0 & 0\\
0 & 0 & 0 & \cdots & \zeta^{n-2} & 0\end{array}\right],\quad\zeta\in\mathbb{T},\]
and $P_{0}$ a diagonal projection. The diagonal elements $(\alpha_{1},\alpha_{2},\dots,\alpha_{n})$
of this projection are defined by setting $\alpha_{i}=1$ if $i\in A$
and $\alpha_{i}=0$ otherwise. The pair $(U_{0}(\zeta),P_{0})$ is
irreducible provided that $n$ is the smallest positive period of
$A$.

\end{document}